\documentclass[11pt]{article} 

 \addtocounter{secnumdepth}{-2}

 \usepackage{amssymb,amstext,amsfonts,amsmath,amscd,amsthm}       
 \usepackage{latexsym}      
 \usepackage{epsfig,epic, color}
 \usepackage{pstricks,pst-node,pst-text,pst-tree}
 \usepackage[all]{xy}
 \usepackage{eucal}

 \newcommand{\numberset}[1]{\ensuremath{\mathbb{#1}}} 
 \newcommand{\C}{\numberset{C}}
 \newcommand{\Z}{\numberset{Z}}
 \newcommand{\R}{\numberset{R}}
 \newtheorem{thm}{Theorem}[section]
 \newtheorem{prop}[thm]{Proposition}
 \newtheorem{lem}[thm]{Lemma}
 \newtheorem{cor}[thm]{Corollary}

 \theoremstyle{definition}
 \newtheorem{rem}[thm]{Remark}
 \newtheorem{ex}[thm]{Example}
 \newtheorem{defi}[thm]{Definition}
 \newtheorem{con}[thm]{Conjecture}
 \newtheorem{propdefi}[thm]{Proposition-Definition}
 \newtheorem{constr}[thm]{Construction}
 \DeclareMathOperator{\im}{Im} \DeclareMathOperator{\re}{Re}

  \DeclareMathOperator{\Arg}{Arg}

\begin{document}
\author{Ricardo Casta\~no Bernard}
\title{Symplectic invariants of some \\ families of Lagrangian $T^3$-fibrations}
\date{}
\maketitle

\abstract{ We construct families of Lagrangian 3-torus fibrations
resembling the topology of some of the singularities in  \textit{Topological Mirror Symmetry} \cite{TMS}. We perform a detailed analysis of the affine structure on the base of these fibrations near their discriminant loci. This permits us to classify the aforementioned families up to fibre preserving symplectomorphism. The kind of degenerations we investigate give rise to a large number of symplectic invariants.}

\section{Introduction}

There is increasing interest in the geometry of Calabi-Yau
manifolds. Much of this interest is motivated by an intriguing
relation between pairs of Calabi-Yau manifolds called Mirror
Symmetry. This relation interchanges-- in a highly non-trivial way
--the complex structure of a Calabi-Yau manifold, $Y$, with the
symplectic structure of its mirror, $\check Y$. There are several
approaches to Mirror Symmetry; one of these is proposed by
Strominger, Yau and Zaslow (SYZ) \cite{SYZ}. The SYZ Conjecture
claims-- based on string theoretic arguments --that the mirror
relation can be explained in terms of certain duality between
$T^n$-fibrations on a pair of mirror Calabi-Yau manifolds.

\medskip
The purpose of this paper is primarily motivated by the SYZ
Conjecture; we are interested in the symplectic geometry of
Calabi-Yau manifolds fibred by tori. Let $(Y,J,\omega)$ be a
compact K\"ahler manifold of $\dim_\C (Y)=n$ with complex
structure $J$ and K\"ahler (symplectic) form $\omega$. We say that
$Y$ is a Calabi-Yau manifold if the canonical line bundle has a
non-vanishing global section $\Omega$ such that
$c\Omega\wedge\bar\Omega=\omega^n$ for some constant $c$
\cite{GHJ}. A very popular example of a Calabi-Yau 3-fold is the
(smooth) quintic hypersurface in $\mathbb{P}^4$ defined by:
\begin{equation}
x_0x_1x_2x_3x_4+ t(x_0^5+x_1^5+x_2^5+x_3^5+x_4^5)=0,
\end{equation}
where $t\in D_0\subset\C$ a small punctured disk around $0$.

\medskip
A submanifold $L$ of $Y$ is called Lagrangian if $\omega |_L=0$
and $\dim_\R L=n$. A Lagrangian submanifold $L$ satisfying
$\im\Omega|_L=0$, is called special Lagrangian. This term was
coined by Harvey and Lawson \cite{HL}.

\medskip
A first attempt to state the SYZ Conjecture using mathematical
language can be outlined as follows (c.f. \cite{G-Wilson1},
\cite{splagI} and \cite{splagII}):

\medskip
\begin{con}\label{SYZ-con} Let $Y$ and $\check Y$ be a mirror pair of Calabi-Yau $n$-folds satisfying certain additional conditions. Let $B$ be a compact connected manifold. Then there exists $C^\infty$ maps $f:Y\rightarrow B$ and $\check f:\check Y\rightarrow B$ such that for $b\in B$ the fibres $f^{-1}(b)$ and $\check f^{-1}(b)$ are special Lagrangian. There is a codimension 2 closed subset $\Delta\subseteq B$ such that the fibres $f^{-1}(b)$ and $\check f^{-1}(b)$ over $b\in B\setminus\Delta$ are dual $n$-tori.
\end{con}

\medskip
The idea of duality in Conjecture \ref{SYZ-con} can be explained
in the following way. First consider the $T^n$-bundle
$f_0:Y_0\rightarrow B_0=B\setminus\Delta$, resulting from removing
the singular fibres of $f$. Let $G$ be a group ($G=\R$ or $\Z$ for
our purposes) and denote by $R^kf_{0\ast} (G)$ the locally
constant sheaf on $B_0$ induced by the presheaf $\mathcal
R^kf_{0\ast}(G) =\{ U\mapsto H^k(f^{-1}(U),G), U\subseteq B\}$.
Denote by $\mathcal E=R^1f_{0\ast} (\R )\otimes C^\infty (B_0)$.
This gives a rank $n$ vector bundle $\mathcal E\rightarrow B_0$
with $ R^1f_{0\ast} (\Z )$ being a family of rank $n$ lattices
lying inside $\mathcal E$. Letting  $\check Y_0=\mathcal E\slash
R^1f_{0\ast} (\Z )$ we can define the dual of $f_0$ as the
$T^n$-bundle:
\begin{equation}
\check f_0:\check Y_0\rightarrow B_0.
\end{equation}
According to Conjecture \ref{SYZ-con} one expects to recover the
mirror of $Y$ as a compactification of $\check Y_0$, obtained by
means of gluing on suitable singular fibres. This method raises a
number of issues demanding careful consideration. As pointed out
in \cite[\S 5]{SYZ}, understanding the structure of the singular
fibres is probably one of the crucial issues.

\medskip
Conjecture \ref{SYZ-con} appears to be the right approach if one
pays attention to the topology only, i.e. forgetting about the
complex and symplectic structures and considers $Y$ and $\check Y$
as $C^\infty$ manifolds only. Under some mild assumptions on the
singular fibres-- they are assumed to be semi-stable, i.e., with
unipotent monodromy  --the SYZ duality explains a topological
version of Mirror Symmetry for the quintic:

\begin{thm}[Gross \cite{TMS}]\label{Thm mark}
Let $Q\subseteq\mathbb{P}^4$ be a smooth quintic 3-fold and $\Xi$
a 4-simplex. There is a $T^3$ fibration $g:Q\rightarrow
\partial\Xi$ with semi-stable fibres. The dual $\check g:\check
Q\rightarrow \partial\Xi$ has only semi-stable fibres and $\check
Q$ is diffeomorphic to a mirror pair of $Q$.
\end{thm}

Both fibrations $g$ and $\check g$ have the same discriminant
locus which consists of a trivalent graph $\Gamma$ lying over the
faces of $\partial\Xi$. There are three types of singular fibres
present in both $Q$ and $\check Q$. Let $s\in\Gamma$ and let $Q_s$
be a singular fibre of either $g$ or $\check g$. Let $(b_1,b_2)$
where $b_i=rankH^i(Q_s,\mathbb{Z})$, $i=1,2$. Then $Q_s$ can be
one of the following types:

\begin{itemize}
\item \textit{type $(2,2)$.} This fibre is $S^1\times I_1$, where $I_1$
is a Kodaira type $I_1$ fibre (a pinched torus).  So, fibres of
type $(2,2)$ are singular along a circle. Fibres of type $(2,2)$
lie over the edges of $\Gamma$;

\item \textit{type $(1,2)$.}  This fibre is obtained by collapsing a torus
$T^2\times\{ p\}$ on $T^2\times S^1$ to a point. Fibres of this
kind lie over some vertices of $\Gamma$;

\item \textit{type $(2,1)$}. Let $S\subset T^2$ be a ``figure eight'' (c.f. \cite[fig. 2.2]{TMS}). This fibre is obtained by collapsing the circles $\{ p\}\times S^1$, $p\in S$, on $T^2\times S^1$ to a point. Fibres of this kind lie over some vertices of $\Gamma$.
\end{itemize}
The $(2,1)$ fibre is dual to the $(1,2)$ fibre (their local monodromy representations are dual), whereas the $(2,2)$ fibre is self-dual.

\medskip
One can try to add on structures to the above topological picture.
First, one can try to put suitable symplectic structures on $Q$
and $\check Q$ making $g$ and $\check g$ into Lagrangian
fibrations. The next step would be to put suitable (almost)
Calabi-Yau structures on $Q$ and $\check Q$ making $g$ and $\check
g$ into special Lagrangian fibrations. Recent development on
special Lagrangian geometry (c.f. Joyce \cite{Joyce-SYZ}) suggests
that this program may not be fully completed in the strong terms
of Conjecture \ref{SYZ-con}.  This is not conclusive, however.

\medskip
There has been some progress in the symplectic category.  Wei-Dong
Ruan \cite{Ruan2}, \cite{Ruan-III} constructs Lagrangian
torus fibrations on the quintic. Ruan's method consists, roughly speaking, on a
certain gradient flow deformation of a well known Lagrangian
fibration on the normal crossing quintic to a neighbour
non-singular quintic.  This produces a piecewise $C^\infty$
fibration with codimension 1 discriminant locus. The topology of this fibration differs from the topological $T^3$ fibration in \cite{TMS}. Ruan argues \cite{Ruan3} that the codimension 1 discriminant can be deformed to codimension 2, in which case the resulting fibration coincides, topologically, with the one in \cite{TMS}.

\medskip
In this paper we are interested in the \textsl{semi-global}
symplectic geometry in a neighbourhood of the singular fibres
rather than in the global picture. We follow the spirit of
\cite{TMS} and construct singular local models of Lagrangian $T^3$
fibrations first. We are able to construct $C^\infty$ Lagrangian
$T^3$ fibrations with singular fibres resembling the topology of
the $(2,2)$ and $(1,2)$ fibres. To date there is no symplectic
model for the $(2,1)$ fibre and it is not yet clear whether there
exists a Lagrangian $T^3$ fibration on the quintic presenting
singular fibres of type $(2,1)$.

\medskip
In two dimensions, it is known that the Kodaira type $I_1$
degeneration of an elliptic fibration, has an infinite number of
symplectic invariants (c.f. V\~u-Ngoc \cite{San} and \cite{RCB}
for an alternative approach and for the $C^k$-symplectic case). We
show that a similar behaviour appears in dimension $n\geq 3$: the
families of Lagrangian $T^n$ fibrations considered in this paper
have infinite dimensional classifying spaces.

\medskip
Gross and Wilson \cite{G-Wilson2} and, independently, Kontsevich
and Soibelman \cite{Kontsevich-Soibelman} propose an alternative
interpretation of SYZ Conjecture, which can be regarded as a
relaxed version of Conjecture \ref{SYZ-con}. This new proposal
posits that the SYZ Conjecture is true in certain limiting sense
as the mirror pair of Calabi-Yau manifolds approach large complex
structure limits. The SYZ duality is then interpreted as a certain
kind of Legendre transform between (singular) affine structures on
the bases of the fibrations. This conjecture is proved for K3
surfaces by Gross and Wilson \cite{G-Wilson2}. For the case of the
quintic, Gross constructs an affine structure on the complement of
the discriminant locus (c.f. \cite[\S19.3]{GHJ}). This affine
structure, in turn, induces a symplectic structure on the
complement of the union of the singular fibres. The results here
can be interpreted as a description of how this affine structure
may become singular at the discriminant locus and the symplectic
invariants arising from this degeneration.

\subsection{Statement of the main results}

\medskip
Let $g:Q\rightarrow \partial\Xi$ be the $T^3$ fibration on the
quintic as in Theorem \ref{Thm mark}. Let $s_0\in \Gamma$ and
consider $U\subseteq \partial\Xi$ an open neighbourhood of $s_0$.
Assume $U$ is small enough so that it contains at most one vertex
of $\Gamma$ and such that $\Gamma\cap U$ is connected. Then
$g^{-1}(U)\subseteq Q$ is a neighbourhood of the fibre
$g^{-1}(s_0)$ and the restriction of $g$ to $g^{-1}(U)$ gives a
$T^3$ fibration, $g^{-1}(U)\rightarrow U$, which is singular along
$\Gamma\cap U$.

\medskip
Now suppose there is a non-compact symplectic manifold $(X,\omega
)$ together with a proper Lagrangian fibration $f:X\rightarrow B$.
In addition, suppose there are diffeomorphisms $\Phi$ and $\phi$
giving a commutative diagram:
\begin{equation}\label{diagram: s-g}
\begin{CD}
X @>\Phi >> g^{-1}(U) \\
@VfVV @VVV \\
B @>\phi >> U
\end{CD}
\end{equation}
Let $\Delta=\phi^{-1}(\Gamma\cap U)$ and
$b_0:=\phi^{-1}(s_0)\in\Delta$. Then $f$ is a Lagrangian $T^3$
fibration having $\Delta$ as discriminant locus and the fibre
$X_{b_0}:=f^{-1}(b_0)$ is homeomorphic to the fibre $g^{-1}(s_0)$.

\begin{figure}[!ht]
\begin{center}
\input{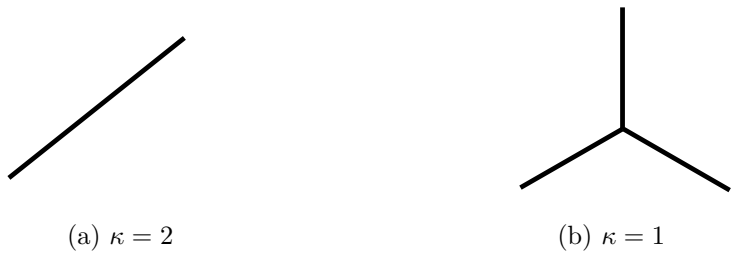}
\caption{The discriminant loci $\Delta\subset B$.}\label{locdis}
\end{center}
\end{figure}

\begin{defi}\label{defi symplecto} Let $\mathcal L(X_{b_0})$ denote the set of triples $\mathcal F=(X,\omega, f)$ where $f:(X,\omega )\rightarrow B$ is a Lagrangian fibration arising as in diagram (\ref{diagram: s-g}). We say that two elements $(X,\omega, f)$ and $(X',\omega', f')$ in $\mathcal L(X_{b_0})$ are \textit{symplectically equivalent} if there is a
symplectomorphism $\Psi :X\rightarrow X'$ and a diffeomorphism $\psi :B\rightarrow B$ such that $\psi (b_0)=b_0$ and $f'\circ\Psi =\psi\circ f$. The set of equivalence classes under
this relation will be denoted by $\widetilde{\mathcal L}
(X_{b_0})$. The elements of $\widetilde{\mathcal L} (X_{b_0})$ can
be regarded as germs of Lagrangian fibrations around the singular
fibre $X_{b_0}$.
\end{defi}
There are three families to be considered: $\mathcal L(2,2)$, $\mathcal L(1,2)$ and $\mathcal L(2,1)$ corresponding to $X_{b_0}$ of type $(2,2)$, $(1,2)$ and $(2,1)$ respectively. Their discriminant loci are as depicted in Figure \ref{locdis}.

\begin{thm}  There is an element $\mathcal F_H\in\mathcal L(\kappa ,2)$, $\kappa=1,2$, for each $H\in C^\infty(B)$.
\end{thm}

\begin{thm}
The germs of fibrations of type $\widetilde{\mathcal L} (\kappa, 2)$, $\kappa=1,2$, are classified by $C^\infty_\Delta (B)$, the space of germs of $C^\infty$ functions on $B$ vanishing at $\Delta$ to all orders.
\end{thm}

\section{Preliminaries}

We recall the construction of action-angle variables on Lagrangian $T^n$ bundles. This is an extensively used technique in the context of Hamiltonian mechanics. The material presented here is standard (c.f. \cite{Dui} and \cite[\S 2]{splagII}).

\medskip
Let $(X,\omega )$ be a symplectic manifold of dimension $2n$ and
$B$ a $n$-dimensional manifold. We shall assume $X$ and $B$ to be
connected but not necessarily compact.

\begin{defi}\label{def uno}
Let $f:X\rightarrow B$ be a proper $C^\infty$ map with connected fibres and denote by $Crit (f)\subset X$ the set of points in $X$ where the differential $f_\ast$ is not surjective. Let $X^\# =X\setminus Crit(f)$ and
$f^\#$ denote the restriction of $f$ to $X^\#$. If the fibres of $f^\#$ are Lagrangian with respect to $\omega$ we say that $f$ is a \textit{Lagrangian fibration}. We denote a Lagrangian fibration as a triple $\mathcal F=(X,\omega ,f)$.
\end{defi}

Observe that the Arnold-Liouville theorem implies that the regular fibres of $\mathcal F$ as in Definition \ref{def uno} are diffeomorphic to $T^n$.

\begin{defi}\label{ass. one} Let $\mathcal F=(X,\omega ,f)$ be a Lagrangian fibration. Denote $X_b=f^{-1}(b)$ and let $Crit(X_b)=Crit(f)\cap f^{-1}(b)$ the set of singular points of $X_b$. We say that $\mathcal F$ is \textit{admissible} if
\begin{itemize}
\item[(1)] $Crit(X_b)$ is connected and the fibres of $f^\#$ are connected;
\item[(2)] $\Delta =f(Crit(f))$ is a closed codimension two subset of $B$ and
\item[(3)] $f^\# (X^\# )=B$ and for any point $x\in X^\#$ there is a local $C^\infty$ section of $f$ passing through $x$.
\end{itemize}
\end{defi}
Observe that (3) in Definition \ref{ass. one} implies that $f$
does not have singular fibres dropping dimension. The fibres of
$f^\#$ over $\Delta$ are diffeomorphic to $T^k\times\R^{n-k}$.
From now on we only consider Lagrangian fibrations which are
admissible.

\medskip
Now let $B_0=B\setminus\Delta$, $X_0=f^{-1}(B_0)$ and $f_0=f|_{X_0}$. The map, $f_0:X_0\rightarrow B_0$ defines a $T^n$ fibre bundle, denoted by $(X_0,f_0)$. Consider $R^{n-1}f_{0\ast}\Z$, a local system as defined in \S 1. Since $f_0$ is proper, one can identify the stalk $(R^{n-1}f_{0\ast}\Z)_b$ with $H_{1}(X_b,\Z )$ using Poincar\'e duality.

\medskip
Now consider $f^\# :X^\#\rightarrow B$ and let $X^\#_b=f^{\#-1}(b)$. We define a sheaf on $B$, with stalk\footnote{Here
$H_c^\ast(\;\cdot\ ,\Z )$ denotes compactly supported cohomology with coefficients in $\Z$} $H^{i}_c(X_b^\#,\Z )$ as follows. Let $U\subseteq B$ and consider the presheaf defined by $U\mapsto H^{i}_c(f^{\#-1}(U),\Z )$.
The latter induces a sheaf, denoted
$R^{i}_cf^\#_\ast \Z$, with stalk $(R^{i}_cf^\#_\ast \Z )_b\cong
H^{i}_c(X_b^\#,\Z )$. Again, we can identify $H_r(X^\#_b,\Z )$
with $(R^{n-r}_cf^\#_\ast \Z )_b$.

\medskip
Now we define a map $R^{n-1}_cf^\#_\ast\Z \hookrightarrow
T^\ast_{B}$ as follows. For each $U\subseteq B$ open and $b\in U$
let $\gamma (b)\in H_1(X_b^\#,\Z )\cong H_c^{n-1}(X_b^\#,\Z )$,
$v\in T_{U,b}$ and $\tilde v$ a lifting of $v$. Define the map
$(b,\gamma (b))\mapsto \lambda_b$, where
\begin{equation}\label{def. sheaf of 1-forms}
\lambda_b (v)= -\int_{\gamma(b)}\iota (\tilde v)\omega.
\end{equation}
This gives a local section $b\mapsto \lambda_b$ of $T_B^\ast$, i.e. a 1-form on $U\subseteq B$. One can check that the above formula does not depend on the lifting of $v$.
\begin{defi} Let $\Lambda \subset T^\ast_{B}$ be the image of $R^{n-1}_cf^\#_\ast\Z$ under the map
 (\ref{def. sheaf of 1-forms}). We call $\Lambda$ the \textit{period lattice}
 of $f$.
\end{defi}

\medskip
Now let us consider $R^{n-1}f_{0\ast}\Z$. Choose a local section $\gamma$ of $R^{n-1}f_{0\ast}\Z$ over an open set on $U\subseteq B_0$. The image of this section under the map
(\ref{def. sheaf of 1-forms}) gives us a period 1-form $\lambda_\gamma$.
This form is closed, since it is the differential of the action
function:
\begin{equation}\label{action int}
\mathcal A_\gamma (b)=\int_{\gamma (b)}\sigma .
\end{equation}
Here $\sigma$ is such that $d\sigma=\omega$. We can ensure that such a $\sigma$ always exists on $f_0^{-1}(U)$ by taking $U\subset B_0$ small enough. This means that the sections of $\Lambda$ are given
locally by the image of a closed 1-form and, in particular, $\Lambda$ is
Lagrangian with respect to the canonical symplectic structure on
$T^\ast_B$.

\medskip
The above construction gives us an exact sequence:

\begin{equation}\label{eq. exact sequence}
\begin{CD}
0@>>> R^{n-1}_cf^\#_{\ast}\Z @>>> T^\ast_{B}
@>>>T^\ast_{B}\slash\Lambda @>>> 0
\end{CD}
\end{equation}

\medskip
\begin{propdefi}\label{def. Jacobian fibration} The sequence (\ref{eq. exact sequence})
defines a symplectic manifold, $J^\#:=T^\ast_B\slash \Lambda$, and
a Lagrangian fibration $\mathcal J_f:J^\#\rightarrow B$ with fibre
$\mathcal J_f^{-1}(b)=T^\ast_{B,b}\slash \Lambda_b$. We call
$\mathcal J_f$ the \textit{Jacobian fibration} of $f$.
\end{propdefi}
\begin{proof} The Lagrangian nature of $\Lambda\subset T^\ast_B$ implies that translations
over $\Lambda$ along the fibres of $T^\ast_B$ are symplectic
transformations. Therefore, $J^\#:=T^\ast_B\slash \Lambda$
inherits the canonical symplectic structure of $T^\ast_B$. The
bundle projection $T^\ast_B\rightarrow B$ induces the map
$\mathcal J_f :J^\#\rightarrow B$. It follows immediately that the
fibres of $\mathcal J_f$ are Lagrangian.
\end{proof}

The following result is deduced from \cite{Dui} (c.f. \cite[\S 2]{splagII}):

\begin{thm}\label{Thm. Jacobian} Let $(X,\omega ,f)$ be a proper Lagrangian fibration. Let $\mathcal J_f:J^\#\rightarrow B$ be the Jacobian fibration of $f$ as defined in (\ref{def. Jacobian fibration}). Then,
\begin{itemize}
\item[(i)] if $f^\#$ has a global section, $\Sigma :B\rightarrow X^\#$, then there is
a fibre preserving diffeomorphism $\Psi :J^\# \rightarrow X^\#$;
\item[(ii)] if $\Sigma$ is Lagrangian, then the diffeomorphism in (i) is a
symplectomorphism.
\end{itemize}
\end{thm}

\medskip
Duistermaat \cite{Dui} observed that a Lagrangian $T^n$ bundle
$f_0:X_0\rightarrow B_0$ has three invariants: its monodromy, its
Chern class and $[\omega ]\in H^2(X_0,\R)$. This tells us that, by
taking $U\subset B_0$ contractible, we can define a set of
action-angle (canonical) coordinates on $f^{-1}(U)\subset X$ which
allows us to write $\omega$ on $f^{-1}(U)$ as the standard
symplectic structure. Furthermore, the action coordinates provide
$B_0$ with an integral affine structure.

\section{The family $\mathcal L(2,2)$}
Lagrangian $T^2$-fibrations with singular fibre of type $I_1$ are
better known in symplectic geometry as \textit{focus-focus}
singularities; they appear in a number of ``physically relevant"
integrable Hamiltonian systems.

\begin{thm}[\cite{RCB}]\label{thm. ff Example}
Let $D\subseteq\C$ be an open disk with coordinates
$s=s_1+\sqrt{-1}s_2$. For any function $h\in C^\infty (D)$ there
is a Lagrangian $T^2$ fibration $\mathcal F=(\bar X,\omega ,f)$
with singular fibre of focus-focus type and whose period lattice
is generated by $\tau_1=-\log |s|ds_1+\Arg (s)ds_2+ dh$ and
$\tau_2=2\pi ds_2$.
\end{thm}

\begin{rem}
There is an alternative proof of Theorem \ref{thm. ff Example}
proposed by V\~u-Ngoc \cite{San}.
\end{rem}

\begin{prop}\label{Prop toy 22}
Let $\bar f: \bar X\rightarrow D$ be a $T^2$ fibration as in
Theorem \ref{thm. ff Example} and $(0,1)$ an open interval. Let
$X=\bar X\times S^1\times (0,1)$ and define $f:X\rightarrow
D\times(0,1)$ to be the composition of the projection onto $\bar
X\times (0,1)$ and $\bar f\times id$. Then, there is a symplectic
structure $\omega$ on $X$ making the fibres of $f$ Lagrangian.
Furthermore the fibres over $\Delta =\{0 \}\times (0,1)$ are
diffeomorphic to $I_1\times S^1$.
\end{prop}
\begin{proof} Let $(r,\theta )$ be coordinates
on $(0,1)\times S^1$. Define $\omega =\bar\omega+dr\wedge d\theta$.
One can verify $f$ is Lagrangian with respect to $\omega$.
\end{proof}

\begin{defi}
Let $M$ be a symplectic manifold and let $f_1,\ldots ,f_n\in C^\infty (M)$ define an integrable Hamiltonian system on $M$. Let $x\in M$ and let $(t_j;x)\mapsto\phi_j^{t_j}(x)$ be the flow generated by the Hamiltonian vector field of $f_j$. We call $(t_1,\ldots ,t_n; x)\mapsto\phi^{t_n}_n\circ\cdots\circ\phi^{t_1}_1(x)$ the \textit{Poisson action} of the system. If all flows $\phi_j$ are complete the Poisson action is an $\R^n$-action on $M$ which preserves the fibres of the map $x\mapsto (f_1(x),\ldots ,f_n(x))$.
\end{defi}

Observe that for $f$ as in Proposition \ref{Prop toy 22} all points $x\in Crit(f)$ are non-degenerate and such that $\textrm{rank}\, f_\ast |_x=1$. Regarding $f$ as an integrable system, one can check that the Poisson orbit of $x$, $\mathcal O_x$, is diffeomorphic to $S^1$ and each point in $\mathcal O_x$ is a rank one critical point. Rank one singular orbits of integrable systems are classified up to fibre preserving symplectomorphism. We state here a special case of a result due to Miranda and Tien-Zung \cite{Zung6}.

\begin{thm}[\cite{Zung6}]\label{ZM}
Let $(M^6,\Omega, h)$ be a (not necessarily proper) Lagrangian
fibration with a non-degenerate rank 1 singular orbit $\mathcal O$
of the Poisson action. Let $D^4$ be a 4-ball and let
$V=D^4\times (-1,1)\times S^1$ be a symplectic 6-manifold with
canonical coordinates $(x_j,y_j,r,\theta )$. There exists a
neighbourhood $U\subseteq M$ of $\mathcal O$ a Lagrangian
fibration, $L:V\rightarrow D\times (-1,1)$, $L(x_j,y_j,r,\theta
)=(q_1(x_j,y_j),q_2(x_j,y_j),r)$, and a fibre preserving
symplectomorphism $\psi :U\rightarrow V$ sending $\mathcal O$ to
$\{x_i=y_i=r=0\}$ and such that $q_i$ can be one of the following
types:
\begin{center}
\begin{tabular}{l l}
elliptic type: & $q_i=x_i^2+y_i^2$\\
\\
hyperbolic type: & $q_i=x_iy_i$\\
\\
focus-focus type: & $\begin{cases}
q_{i}=x_iy_i+x_{i+1}y_{i+1} \\
q_{i+1}=x_iy_{i+1}-x_{i+1}y_i
\end{cases}$
\end{tabular}
\end{center}
\end{thm}

\begin{defi}\label{defi one} Let $\mathcal F= (X,\omega ,f)$ be an admissible  Lagrangian $T^3$ fibration. Let $x\in Crit(X)$ be a non-degenerate rank 1 singular point and $\mathcal O_x$ its Poisson orbit. We say that $\mathcal F$ is a \textit{Lagrangian fibration of type $(2,2)$}, denoted $\mathcal F\in\mathcal L(2,2)$, if there is a neighbourhood $U\subset X$ of $\mathcal O_x$ such that $X=f^{-1}(f(U))$ and the following commutative diagram:
\begin{equation}\label{eq. 22 normfor}
\xymatrix{
  U \ar[d]_{f|_U} \ar[rr]^{\psi}\ar[drr]^F & & V \ar[d]^{q=(q_1,q_2,q_3)} \\
  f(U) \ar[rr]^{\phi} & & \ \, D\times (0,1)   }
\end{equation}
where $\psi :U\rightarrow V=D^4\times (0,1)\times S^1$ is a symplectomorphism, $\phi$ is a diffeomorphism and $q_1$, $q_2$ are of focus-focus type, $q_3=r$.
\end{defi}

We shall denote $B=D\times (0,1)$ and $b=(b_1,b_2,b_3)\in B$. We can write $\phi =(\phi_1,\phi_2,\phi_3)$ and $\phi\circ f=(f_1,f_2,f_3)$, where $f_i=\phi_i\circ f$. If we think of $\psi$ as providing $U$ with canonical coordinates, then $f_j|_U=q_j$ or, with slight abuse of notation, $f|_U=F$ where $F(x_i,y_i)=(q_1,q_2,q_3)$ as in (\ref{eq. 22 normfor}). We regard $F$ as the normal form for the family $\mathcal L(2,2)$.

\medskip
Let $v_{q_j}$ be the Hamiltonian vector field corresponding to $q_j$ and let $g_j^t$ its flow. Let $\zeta_1=x_1+\sqrt{-1}x_2$ and $\zeta_2=y_1+\sqrt{-1}y_2$. Observe that $\bar\zeta_1\zeta_2=q_1+\sqrt{-1}q_2$. The flows of $g_j:\R\times V\rightarrow V$ are given by:
\begin{equation}\label{eq. 22 flows}
\begin{tabular}{l}
$g_1^t(\zeta_1,\zeta_2,r,\theta )=(e^{t}\zeta_1,e^{-t}\zeta_2, r,\theta)$  \\
$g_2^t(\zeta_1,\zeta_2,r,\theta )=(e^{it}\zeta_1,e^{it}\zeta_2,r,\theta )$\\
$g_3^t(\zeta_1,\zeta_2,r,\theta )=(\zeta_1,\zeta_2,r,\theta -t)$.
\end{tabular}
\end{equation}
Observe that $g_2^t$ and $g_3^t$ generate a fibre-preserving $T^2$ action on $V$.

\begin{lem}\label{lem triviality} Let $(X,\omega ,f)\in\mathcal L(2,2)$. Then the compact fibres of $f^\#:X^\#\rightarrow B$ are diffeomorphic to $T^3$ whereas the non-compact ones are diffeomorphic to $T^2\times\R$. There is an open neighbourhood $\mathcal U\subset X$ of $Crit (f)$ such that the fibres of $f_{\mathcal U}:=f|_{X\setminus\mathcal U}$ are diffeomorphic to $T^2\times [0,1]$. Furthermore, $f_{\mathcal U}$ defines a trivial fibre bundle.
\end{lem}
\begin{proof} The first part follows directly from the definition. For the second claim it is enough to take $\mathcal U$ a small connected neighbourhood of $Crit(f)$ which is invariant with respect to the $T^2$-action induced by $v_{q_2}$ and $v_{q_3}$ and redefine $B:=\phi\circ f(\mathcal U)$. The triviality of $f_{\mathcal U}$ follows from the fact that $B$ is contractible.
\end{proof}

Notice that $f_j|_{\mathcal U}=q_j$ implies that the vector fields $v_{q_j}$ extend vector Hamiltonian fields $v_j$ on $X$ which are tangent to the fibres, hence the flows $g_i^t$ extend to $X$. Since the fibres of $f$ are compact $g_i^t$ are complete.

\begin{constr}\label{contr. action} Define an action $\Pi :\R^3\times X\rightarrow X$, $(T,x)\mapsto\Pi^T(x)$, $T=(t_1,t_2,t_3)$ as the composition of flows:
\[
\Pi (T,x):=g_3^{t_3}\circ g_2^{t_2}\circ
g_1^{t_1}(x).
\]
The restriction of $\Pi$ to $X^\#$ is just the Poisson action on
$f^\# :X^\#\rightarrow B$, which is free and transitive along the fibres
since $f^{\# -1}(b)$ is connected. This implies that for any two
points $x,y\in f^{\# -1}(b)$, there is a multi-time
$T=T(x,y)\in\R^3$ such that $\Pi^{T}(x)=y$.
Similarly, consider now the Hamiltonian vector fields $v_{q_1},v_{q_2},
v_{q_3}$, on $\mathcal U\subseteq X$. In an analogous way, we can define an action on $\mathcal U$, $\Pi_0:R\times \mathcal U\rightarrow\mathcal U$, where $R\subseteq\R^3$ is some open set, as the composition of flows of $v_{q_i}$. Since $F^{-1}(b)$ is connected and non-singular for $b\in B_0=B\setminus\Delta$, $\Pi_0$ is transitive along the regular fibres of $F$.
\end{constr}

We are going to use the actions $\Pi$ and $\Pi_0$ to compute the period lattice of $(X,\omega ,f)$. Let $\epsilon >0$ and write $s=b_1+\sqrt{-1}b_2$, $r=b_3$. Define,
$\Sigma_1(s,r)=(\bar s/\epsilon ,\epsilon ,r,\theta_0)$ and $\Sigma_2 (s,r)=(\epsilon,s/\epsilon ,r,\theta_0)$,
$\theta_0\in S^1=\R\slash\mathbb{Z}$. These give sections of $F$
which lie inside $\mathcal U$ and do not intersect $Crit(F)$. Now consider the equation:
\begin{equation}\label{period-equation 2}
\Pi_0(T_0(b),\Sigma_1(b))=\Sigma_2(b), b\in B_0.
\end{equation}
The solution $T_0=(\alpha_1,\alpha_2, \alpha_3)$ is determined by the system:
\begin{align}
\begin{cases}
e^{-\alpha_1+i\alpha_2}\cdot\epsilon =s/\epsilon \notag \\
e^{\alpha_1+i\alpha_2}\cdot\bar s/\epsilon =\epsilon \notag \\
\theta_0 -\alpha_3=\theta_0  \notag
\end{cases}
\end{align}
One verifies that the (primitive) solution to the system is
\[
\alpha_1=-\log |s|+2\log\epsilon ,\quad \alpha_2=\Arg (s),\quad
\alpha_3=0. \]

\medskip
Let $\mathcal U'$ be a $T^2$ invariant neighbourhood of $Crit(f)$,  $\mathcal U'\subset\mathcal U$ as in Lemma \ref{lem triviality}. We can take $\mathcal U'$ small enough so that we can regard $\Sigma_1$ and $\Sigma_2$ also as sections of the $T^2\times\R$ fibre bundle $f_{\mathcal U'}$ over $B$.

\begin{prop}\label{uniqueness of alpha} Let $\Sigma_1$ and
$\Sigma_2$ be sections of $f_{\mathcal U'}$ as above. The equation:
\begin{equation}\label{period-equation}
\Pi^{T(b)}(\Sigma_2(b))=\Sigma_1(b),\qquad b\in B
\end{equation}
has a unique solution, $T(b)=(\eta_1(b),\eta_2(b),\eta_3(b))$,  which depends smoothly on $b\in B$.
\end{prop}
\begin{proof}
A solution to equation (\ref{period-equation}) exists since the
action $\Pi$ is transitive along the fibres of $f_{\mathcal U'}$. We
shall see that $T(b)$ depends smoothly on $b\in B$.

\medskip
Let $S_j=\Sigma_j(B)$ and let $x_0\in S_2$ such that
$f(x_0)=b_0\in B$. Then, there is $t_0\in\R^3$ such that
$\Pi^{t_0}(x_0)=y_0\in S_1$. Let $U_0$ be a small neighbourhood
of $x_0$ and let $R$ be a neighbourhood of $t_0$. Let $V_0$ be a
neighbourhood of $y_0$ such that $f(V_0)=f(U_0)\subseteq B$.
Define $P:R\times (U_0\cap S_2)\rightarrow V_0$, as
$P(t,x)=\Pi^t(x)$. Notice that $S_1$ is transversal to the
$\R^3$-orbit of $\Pi$ passing through a point $y\in S_1$. This
implies that $P$ is transversal to $S_1\cap V_0$. Then, $\mathcal
P:=P^{-1}(S_1\cap V_0)$ is a codimension $3$ smooth submanifold of
$R\times (U_0\cap S_2)$.

\medskip
Now observe that since $\Pi$ is an
action, the ``time'' derivative of $P$ evaluated at $(t_0,x_0)$ is
non-singular. Then, $\mathcal P$ can be described locally as the
graph of a $C^\infty$ map, $g:U'_0\rightarrow R$, where
$U'_0\subseteq (S_2\cap U_0)$ is a small neighbourhood of $x_0$.
Let $B'=f(U_0')$ and define $T:B'\subseteq B\rightarrow
R$, as $T(b)=g\circ \Sigma_2(b)$ for $b\in B'$. Then, $T$ is a $\R^3$-valued $C^\infty$ function such that
$\Pi^{T(b)}(\Sigma_2 (b))=\Sigma_1(b)$.
From Lemma \ref{lem triviality} we know that $f_{\mathcal U'}$ has trivial monodromy. Therefore these local solutions can be glued together to give a single-valued global solution, $T(b)=(\eta_1 (b),\eta_2(b),\eta_3(b))$, $b\in B$.
\end{proof}

Define the 1-form $\eta =\eta_1db_1+\eta_2db_2+\eta_3db_3$ on $B$ where $\eta_i\in C^\infty (B)$ are as in Proposition \ref{uniqueness of alpha}.

\begin{prop}\label{prop. c1 period lattice} Let $\mathcal F=(X,\omega ,f)\in\mathcal L(2,2)$. There are local sections
$(e_1,e_2,e_3)$ of $R^2_cf^\#\Z$ such that the period lattice of $\mathcal F$ is generated by the 1-forms:
\begin{align}
\tau_1 =\tau_0 +dH,\quad \tau_2 = 2\pi ds_2,\quad  \tau_3 = dr
\notag
\end{align}
where $\tau_0=-\log |s|ds_1+\Arg (s)ds_2$ and $H$ is a smooth
function of $b=(s_1,s_2,r)\in D\times I$ such that $dH=\eta$. The
monodromy of $f$ expressed in terms of $\Lambda
=\langle\tau_1,\tau_2,\tau_3\rangle$ is represented by the matrix:
{\small $\begin{pmatrix}
  1 & 0 & 0\\
  1 & 1 & 0\\
  0 & 0 & 1
\end{pmatrix}$.}
\end{prop}
\begin{proof} One can construct generators
$e_1(b),e_2(b),e_3(b)$ of $H_1(X_b,\Z )$, $b\in B_0$ by
means joining integral curves of $v_j$ in a suitable way. For
instance, we define a representative of $e_1$ to be the ordered
composition of paths $\gamma= (\gamma_1, \gamma_2,\gamma_3,\tilde\gamma_1, \tilde\gamma_2, \tilde\gamma_3)$.
Here, $\gamma_i$ is an integral curve of $v_{q_i}$ starting at a
point $x_{i-1}$ running a time $t_i\in [0, \alpha_i]$ and
finishing at a point $x_i$. Similarly, $\tilde\gamma_i$ is an
integral curve of $v_i$ starting at a point $\tilde x_{i-1}$
running a time $\tilde t_i\in[0, \eta_i]$ and finishing at a point
$\tilde x_i$. Then, the curve $\gamma$ is determined by the
initial condition $x_0=\Sigma_1(b)$, $\tilde x_0=x_3$. It follows
from equations (\ref{period-equation
 2}) and (\ref{period-equation}) that $\gamma$ is closed and non-trivial. For constructing a
representative of $e_j$, $j=2,3$, we take an integral curve
of $v_j$ starting at $\Sigma_1(b)$ and flowing from time $0$ to
$1$. Now we can use formula (\ref{def. sheaf of 1-forms}) to compute the period 1-forms. It follows that
$\tau_1=\sum\alpha_jdb_j +\sum\eta_jdb_j$. Since $\tau_1$ and $\tau_0=-\log |s|ds_1+\Arg (s)ds_2$ are closed, then $\tau_1=\tau_0+dH$ for some $H\in C^\infty (B)$. The computation of $\tau_2$ and $\tau_3$ is direct from (\ref{def. sheaf of 1-forms}).
\end{proof}

\section{The family $\mathcal L(1,2)$}

There is a fairly complete understanding of the class of non-degenerate (Morse-Bott) singularities of integrable Hamiltonian systems  (cf. Eliasson \cite{Eliasson}, Tien-Zung and Miranda \cite{Zung6}). Generically, the function components of the fibration-- i.e. the integrals of the system --can be reduced to quadratic polynomials. In contrast, for some special Lagrangian singularities arising from integrable Hamiltonian systems, one should expect cubic terms (c.f. Fu \cite{Fu}).

\subsection{$T^2$-symmetric special Lagrangian singularities}

\medskip
Let $X$ be a symplectic 6-manifold and $f:(X,\omega )\rightarrow
B$ a Lagrangian fibration which is admissible in the sense of
Definition \ref{ass. one}. Denote by $\omega_0=\sum_idx_i\wedge
dy_i$ the standard symplectic form on $\C^3\cong\R^6$ with
canonical coordinates $(x_i,y_i)$, $z_i=x_i+\sqrt{-1}y_i$ and let
$\Omega_0=dz_1\wedge dz_2\wedge dz_3$.

\begin{defi} Let $f:(X,\omega )\rightarrow B$ be a Lagrangian fibration and let $p\in Crit(f)$ and let $k=rank f_\ast |_p$. Let $O_p$ denote the Poisson orbit of $p$. We say that $O_p$ is a \textit{rank $k$ complexity one singularity} if there is an open neighbourhood $W\subseteq X$ of $O_p$ and a Hamiltonian $T^2$ action $\Phi:T^2\times W\rightarrow W$ such that $f(\Phi(t, x))=f(x)$ for each $(t,x)\in T^2\times W$.
\end{defi}
\begin{rem}
Regarding $f|_W$ as an integrable Hamiltonian system, we can think of $T^2$ as a symmetry group of the system. Notice that if $k=0$, i.e. $p$ is a fixed point of the Poisson action, then $p$ is also a fixed point of $\Phi$. It is a standard fact that the Hamiltonian action of a $k$-torus on a symplectic manifold $M$ is completely determined on a neighbourhood of a fixed point $x_0$ by the \textit{weights} of the isotropy representation of the linear action of $T^k$ on $T_{x_0}M$. These are elements  $\rho_1(x_0),\ldots\rho_n(x_0)\in\mathfrak t^\ast=Lie(T^k)^\ast$ (c.f. Guillemin and Sternberg \cite{Gui-Stern}).
\end{rem}

\begin{defi} Let $f:(X,\omega )\rightarrow B$ be a Lagrangian fibration, $p\in Crit (f)$ and $O_p\in X$ a rank $k$ Poisson orbit of $f$. We say that $O_p$ is \textit{special} if there is an open neighbourhood $U\subseteq X$ of $O_p$ and canonical coordinates $(z_i,\bar z_i)$ on $U$ such that $f|_U$ is special Lagrangian with respect to $(\omega_0,\Omega_0)$.
\end{defi}
\begin{ex} Let $(X,\omega, f)\in\mathcal L(2,2)$ and $p\in Crit(f)$. Then  $O_p\subset X$ is a special singularity, which is also a rank one complexity one singularity.
\end{ex}

\subsection{The Harvey-Lawson singularity}
We review an example proposed by Harvey and Lawson (c.f. \cite[\S III.3.A]{HL}). This provides an example of a special rank zero complexity one singularity.

\medskip
Let us consider the map $F:\C^n\rightarrow\R^n$ given by
$F=(F_1,\ldots ,F_n)$ where
\begin{equation}\label{c1 H-Lawson}
F_1=\im\prod z_i,\ F_k=|z_1|^2-|z_k|^2,\quad  k=2,\ldots n.
\end{equation}
The fibres of $F$ are Lagrangian with respect to the standard symplectic form on $\C^n$; for this, one only needs to check that $\{ F_i,F_j\}=0$ for
$i,j=1,\ldots ,n$. In other words, $F$ defines an integrable
Hamiltonian system. One can also check that $\re (\det_\C (\partial_{\bar z_j}F_i))=0$, hence the fibres of $F$ are special Lagrangian. We observe that the map $\mu:=(F_2,\ldots ,F_n)$ is the moment map of the Hamiltonian
$T^{n-1}$ action on $\C^n$ given by:
\[
(z_1,\dots ,z_n)\mapsto (e^{i\theta_1}z_1,\ldots
,e^{i\theta_n}z_n),
\]
with $\theta_1+\cdots +\theta_n=0$. This action preserves the
fibres of $F$. Now let $x=(x_1,\ldots ,x_n)\in\R^n$ and let $z\in
F^{-1}(x)$. Denote by $T\cdot z$ the $T^{n-1}$-orbit of $z$. Then,
$T\cdot z$ is homeomorphic to $T^{n-1}$ unless $z\in Crit
(F)=\bigcup_{1 \leq i<j\leq n}P_{ij}$ where,
\begin{equation}\label{c1 H-L discriminant}
P_{ij}=\{ (z_1,\ldots ,z_n)\in \C^n\ | \ z_i=z_j=0 \}.
\end{equation}
For $z\in Crit(F)$, the orbit $T\cdot z$ is a torus of lower
dimension and it is a point when $z=0$. A fibre $F^{-1}(x)$
disjoint from $Crit(F)$ is homeomorphic to $T^{n-1}\times\R$ and
for $x\in\Delta := F(Crit(F))$ the fibre $F^{-1}(x)$ is a singular
fibre.

\medskip

\begin{figure}[!ht]
\begin{center}
\input{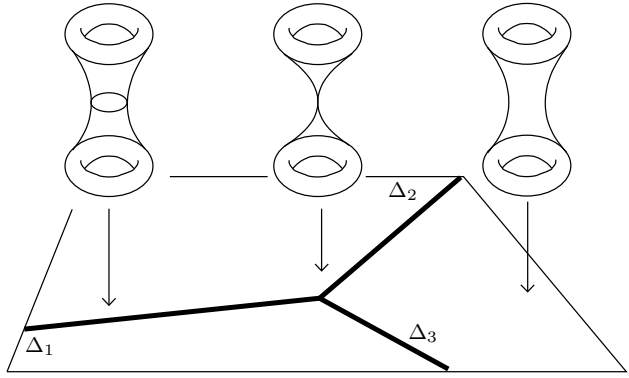}
\caption{The fibres of $F$}\label{Figure HL}
\end{center}
\end{figure}

For $n=3$, $\Delta =\Delta_1\cup\Delta_2\cup\Delta_3\cup\{ 0\}$
where $\Delta_1=\{x_1=0, x_2=x_3> 0\}$, $\Delta_j=\{x_1=x_j=0,
x_{j}<0\}$, for $j=2,3$. The fibre over $x\in\Delta_i$ is
homeomorphic to
\[
S^1\times [\R\times S^1\slash (\{ point\}\times S^1) ],
\]
whereas the fibre over $0\in\R^3$ is homeomorphic to
\[
\R\times T^2\slash (\{ point\}\times T^2)
\]
In particular, we conclude that the map $F$ is not proper.

\begin{rem}\label{rem joyce} Joyce observed (c.f. \cite[\S 5]{Joyce2} and \cite[\S 4]{Joyce-SYZ}) that, in three dimensions, any connected special Lagrangian 3-manifold in $\C^3$ which is invariant under the above $T^2$-action is a subset of some fibre of the map (\ref{c1 H-Lawson}).
\end{rem}

\subsection{The topological $(1,2)$ fibre}

We outline Gross' construction of a topological 3-torus fibration
with fibre of type $(1,2)$. For the details we refer the reader to
\cite[Example 2.10]{TMS}.

\begin{constr}[Gross \cite{TMS}]\label{c1 const 1,2}
Let $B=B^3$ be a 3-ball. We define $\Delta \subset B$ a cone over
three points as follows. Identify $B\setminus\{0\}$ with
$S^2\times (0,1)$ and let $p_1,p_2,p_3\in S^2$. Define
$\Delta_i=\{p_i\}\times (0,1)$. These are the ``legs'' of the
cone. Define $\Delta =\Delta_1\cup\Delta_2\cup\Delta_3\cup\{ 0\}$,
where $\{ 0\}$ is the vertex of the cone.

\medskip
Let $Y=S^1\times B$ and $Y'=Y\setminus (\{
p\}\times\Delta )$, where $p\in S^1$. Let $L\cong\Z^2$ and define $T(L)=L\otimes_{\Z}\R\slash L$. Now consider a principal
$T(L)$-bundle $\pi' :X'\rightarrow Y'$ with Chern class $c_1\in H^2(Y',L)$.
Then the class $c_1$ is represented by a triple $(a_1,a_2,a_3)$ where $a_i\in L$. It is shown (c.f. \cite[Ex. 2.10]{TMS}) that by choosing $c_1=((1,0),(0,1),(-1,-1))$ there is a unique manifold $X$ such that $X'\subset X$ and a commutative diagram of smooth maps:
\[
\xymatrix{
  X' \ar[d]^{\pi'} \ar@{^{(}->}[r] & X \ar[d]^{\pi} \\
  Y' \ar@{^{(}->}[r] & Y   }
\]
such that $\pi$ is proper. Furthermore, it is shown that in a
neighbourhood $U\cong\C\times\R^2\subset Y$ of the vertex of
$\Delta$, the map $\pi:\pi^{-1}(U)\rightarrow U$ coincides
with the map $\tilde\pi :\C^3\rightarrow\C\times\R^2$ given by:
\begin{equation}
\tilde\pi (z_1,z_2,z_3 )=\left (
z_1z_2z_3,|z_1|^2-|z_2|^2,|z_1|^2-|z_3|^2\right ).
\end{equation}

Now define $f:X\rightarrow B$ to be the composition of
$\pi :X\rightarrow Y$ with the projection $Y\rightarrow B$. Then, $f$ is a continuous map whose fibre over
$b\in B\setminus\Delta$ is $T^3$. The fibre over $b\in\Delta_i$ is
homeomorphic to $S^1\times [S^1\times S^1\slash (\{ point\}\times S^1) ]$,
i.e. it is a $(2,2)$ fibre, whereas the fibre over the vertex of
$\Delta$ is homeomorphic to $S^1\times T^2\slash (\{ point\}\times T^2)$,
i.e. it is a $(1,2)$ fibre.
\end{constr}

\begin{figure}[!ht]
\begin{center}
\input{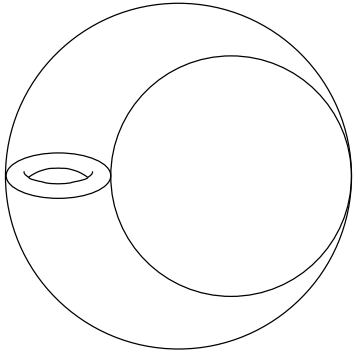}
\caption{Singular fibre of type $(1,2)$}\label{one-two}
\end{center}
\end{figure}

It turns out that the $T(L)$ action on $X'$ action extends to $X$, moreover, $Crit(f)\subset X$ consists of the union of the critical orbits of this action. There is a single fixed point $p\in Crit(f)$, which is singular point of $f^{-1}(0)$.

\subsection{The symplectic structure}
Let $f:X\rightarrow B$ as in Construction \ref{c1 const 1,2} and suppose there is a symplectic structure $\omega$ making $f$ Lagrangian, i.e. defining a triple $(X,\omega, f)\in\mathcal L(1,2)$. Furthermore, assume the extended $T(L)$ action on $X$ preserves $\omega$. It is follows from these hypotheses that $p$ is a rank zero complexity one singularity. Let $\mathfrak t=Lie(T(L))$. Then we can regard $L\hookrightarrow\mathfrak t$ and identify $c_1$ with the isotropy data $(\rho_1(p),\rho_2(p),\rho_3(p))$ of the $T(L)$-action at $p$.

\begin{thm}\label{thm. norm-form}
Let $f:(X,\omega)\rightarrow B$ be a Lagrangian fibration of type
$(1,2)$. Assume there is a fibre-preserving $T:=T(L)$ action
preserving $\omega$. Let $p\in Crit (f)\cap f^{-1}(0)$. Then there
is an open  neighbourhood $U\subset X$ of $p$, a symplectomorphism
$\psi :U\rightarrow V\subseteq (T_pX,\omega_0)$ and a
diffeomorphism $\varphi :f(U)\rightarrow \R^3$ such that
$q\circ\psi=\varphi\circ f|_U$, where $q:V\rightarrow\R^3$ is a
Lagrangian fibration given by
$q=(h,|z_1|^2-|z_2|^2,|z_1|^2-|z_3|^2)$, $h\in C^\infty (V)$.
Furthermore, if $p$ is special, then $h=\im z_1z_2z_3$.
\end{thm}
\begin{proof}
Consider $\mu :X\rightarrow\mathfrak{t}^\ast$ the moment map of
the $T$ action around $p$. According to \cite{Gui-Stern}, there is
a neighbourhood $U\subseteq X$ of $p$ and an equivariant
symplectomorphism $\psi :U\rightarrow V\subseteq (T_pX,\omega_0)$
such that $\mu=\psi^\ast M$, where $M =c+\sum_i \rho_i|z_i|^2$,
$c\in\mathfrak t^\ast$. Without loss of generality we can assume
$c=0$ and choose a basis of $\mathfrak t^\ast$ such that
$\rho_1=(1,1)$, $\rho_2=(-1,0)$ and $\rho_3=(0,-1)$. Then we can
write $M=(M_1,M_2)$, where $M_j=|z_1|^2-|z_j|^2$. Let $v_j$ be the
vector field on $V$ determined by the equation: $\iota
(v_j)\omega_0=dM_j$. The orbits of $v_j$ are periodic with period
$2\pi$. Now let $\Sigma$ be a section of $f$ over $B':=f(U)$ such
that $\Sigma (B')\subset U\setminus Crit (f)$. Let $y(b)=\psi
(\Sigma (b))$ and $g_j:[0,2\pi ]\rightarrow V$ be an integral
orbit of $v_j$ passing through $y(b)$. Then $g_j$ pulls back to a
loop $\gamma_j(b)\subset f^{-1}(b)\cap U$, disjoint from
$Crit(f)$. We can assume there is a 1-form $\sigma$ such that
$d\sigma =-\omega$. Let $\mathcal
A_j(b)=\int_{\gamma_j(b)}\sigma$. One can verify that $\mathcal
A_j\circ f|_{U}=M_j\circ\psi$. Now let $\alpha$ be a Lagrangian
section of $T^\ast_{B'}$ close to the zero section. We can choose
$\alpha$ such that $\alpha (0)\wedge d\mathcal A_1(0)\wedge
d\mathcal A_2(0)\neq 0$. Then there is an open neighbourhood of
$0\in B'$ in which $\alpha\wedge d\mathcal A_1\wedge d\mathcal
A_2\neq 0$ and a unique smooth function $\mathcal A$ such that
$\mathcal A (0)=0$ and $d\mathcal A =\alpha$. Then $\varphi
=(\mathcal A,\mathcal A_1, \mathcal A_2)$ defines a
diffeomorphism from a small neighbourhood of $0\in\R^3$ denoted,
with abuse of notation by $B$, into $\R^3$,
$\varphi:B\rightarrow\varphi (B)\subseteq\R^3$. Let $h=\mathcal
A\circ f\circ\psi^{-1}$. Then $h$ is a $T$-invariant function on
$V$ hence $q:=(h,M_1,M_2)$ defines a Lagrangian fibration on $V$
such that $q=\varphi\circ f|_U\circ\psi^{-1}$. Now we can think of
$\psi$ as identifying $U\cong V\subseteq\C^3$ such that
$\varphi\circ f|_U=(h,M_1,M_2)$. In view of Remark \ref{rem joyce}
and since $\psi (Crit(f))=Crit(q)=Crit (M)=\bigcup_{ij}\{
z_i=z_j=0\}$,  we see that if $\varphi\circ f|_U$ is special
Lagrangian then there should exist a 1-form $\alpha$ with the
above properties and such that $h=\im z_1z_2z_3$.
\end{proof}

\begin{rem}
Observe that the $T(L)$-action on $X$ can always be assumed to be
Hamiltonian. Indeed, the above action is chosen so that $f$ has
the desired monodromy, in particular, it induces monodromy
invariant cycles $e_1(b),e_2(b)\in H_1(f^{-1}(b),\Z)$ which can be
used to compute the action integrals $\mathcal A_{e_1}$, $\mathcal
A_{e_2}$. Then $\mu_i=\mathcal A_i\circ f$ define the moment map
$(\mu_1,\mu_2)$ of a $T=S^1\times S^1$ action, which is defined on
$X$ as $e_i$ are monodromy invariant. It is a consequence of
\cite[Prop. 3.3]{TMS}  and \cite[Thm. 2.2]{Gross_spLagEx} that
$p\in f^{-1}(0)\cap Crit(f)$ can be made into a special
singularity with respect to $(\omega_0,\Omega_0)$.
\end{rem}

\begin{cor}\label{cor norm form 1,2} Let $(X,\omega, f)\in\mathcal L(1,2)$. Let $p\in f^{-1}(0)\cap Crit(f)$ and let $(T_pX,\omega_0)$, where $\omega_0=\frac{i}{2}\sum dz_j\wedge d\bar z_j$. There is a neighbourhood $U\subset X$ of $p$ and a 3-ball $B$ centred at $0\in\R^3$ such that $f(U)=B$, a diffeomorphism $\varphi :B\rightarrow \varphi (B)\subseteq\R^3$ and a symplectomorphism $\psi :U\rightarrow V\subseteq (T_pX,\omega_0)$ such that $\varphi\circ f|_U=F\circ\psi$ where $F(z_1,z_2,z_3)=(\im z_1z_2z_3,|z_1|^2-|z_2|^2,|z_1|^2-|z_3|^2)$.
\end{cor}

\medskip
\subsection{Example of a Lagrangian fibration of type $(1,2)$}

Here we show that the family $\mathcal L(1,2)$ is not void. We construct a Lagrangian fibration of type $(1,2)$ for each $H\in C^\infty(B)$, $B\subset\R^3$ an open ball. The arguments we use here are valid in any dimension.

\medskip
Consider the map $F:\C^n\rightarrow\R^n$, where $F=(F_1,\ldots
,F_n)$ as in (\ref{c1 H-Lawson}). The quotient $\C^n\slash T^{n-1}$ can be identified with $\C\times\R^{n-1}$ by means of the map $\pi :\C^n\rightarrow\C\times\R^{n-1}$,
\begin{equation}\label{c1 pi on V}
\pi (z)=(\prod z_i,|z_1|^2-|z_2|^2,\ldots ,|z_1|^2-|z_n|^2).
\end{equation}
Let $\prod z_i=u+\sqrt{-1}b_1\in\C$ and $b_j=|z_1|^2-|z_j|^2$ and $b=(b_1,\ldots b_n)$. Letting $x_i=|z_i|^2$ the following relations hold:

\begin{equation}{\label{circle-bundle-eq}}
  \begin{cases}
     & \prod_{i=1}^nx_i=u^2+b_1^2, \\
     & x_1-x_j=b_j,\quad j\geq 2 . \\
  \end{cases}
\end{equation}
We can restate these equations (renaming $x:=x_1$) as:
\begin{equation}{\label{polynomial}}
x\prod_{j\geq 2}(x-b_j)-b_1^2=u^2.
\end{equation}
Define $P_b(x)=x\prod_{j\geq 2}(x-b_j)-b_1^2$. We can regard
$P_b(x)$ as a polynomial in the variable $x$ with $b\in\R^n$
acting as a parameter. We notice that for all values of $b$,
$P_b(x)=0$ has always a non-negative \textit{real} solution.
Define $\mathcal Z^\R_b=\{ \zeta(b)\in\R \mid P_b(\zeta (b))=0 \}$.
This is an ordered set, so we can take $\zeta_0(b)=\max\mathcal
Z^\R_b$. We observe that $P_b(x)>0$ for $x>\zeta_0(b)$; $P'_b(x)\neq 0$ for $x>\zeta_0(b)$ and $P_b'(\zeta_0(b))=0$ if and only if $b\in\Delta$. Observe that $\zeta_0(b)$ becomes a multiple root of $P_b(x)$ when $b\in\Delta$.

\begin{lem}\label{estimates of zeta} The function $\zeta_0(b)$ is smooth on $\R^n\setminus\Delta$ and continuous on $\Delta$. Let
 $\partial^k_{J_k}\zeta_0$ denote an order $k$ partial derivative of $\zeta_0$, $J_k=j_1,\ldots ,j_n$, $j_1+\cdots+j_n=k$. Let $B\subset\R^n$ be a small neighbourhood of $0\in\R^n$. Then,
\begin{equation}\label{partials of zeta}
\partial^k_{J_k}\zeta_0=\sum_{l<\infty}\frac{G_l(b,x)|_{\zeta_0}}{(P'_b(\zeta_0))^{\lambda_l}},
\end{equation}
where $G_l(b,x)$ is bounded on $B$ and $\lambda_l\in\Z^+$ is a finite
power.
\end{lem}
\begin{proof}
For $b\in\R^n\setminus\Delta$, $P'_b(\zeta_0(b))\neq 0$ and it follows that $\zeta_0(b)$ is smooth on $\R^n\setminus\Delta$. Let $G(b)=P_b(\zeta_0(b))$ and consider
$\partial_{b_j}G$. We notice that $G\equiv 0$ on $\R^n$, hence
$\partial_{b_j}G\equiv 0$. This implies
\[
\partial_{b_j}\zeta_0=-\frac{\partial_{b_j}P_b\big|_{\zeta_0}}{P'_b(\zeta_0)}.
\]
The function $G_{1_j}(b)=\partial_{b_j}P_b\big|_{\zeta_0}$ is bounded on a small ball $B\subset\R^n$ centred at $0\in\R^n$. The verification of the case $k>1$ is left to the reader.
\end{proof}
Let  $\epsilon >0$ and let $\zeta_\epsilon (b)$ be the maximal real
solution of $P_b(x)-\epsilon^2=0$. Observe that $\zeta_0(b)<\zeta_\epsilon(b)$ and $P'_b(\zeta_\epsilon(b))\neq 0$ for all $b\in\R^n$. It is easy to verify that $\zeta_\epsilon(b)$ is a smooth function on $\R^n$.
\begin{cor}\label{cor sect HL} Let $F:\C^n\rightarrow \R^3$ as in (\ref{c1 H-Lawson}). Let $\zeta_0$ and $\zeta_1$ be the maximal real solutions of $P_b(x)=0$ and $P_b(x)-1=0$ respectively. Let $\theta_\pm (b)=\Arg (\pm 1 +ib_1)$. The maps $\Sigma^-$ and $\Sigma^+$,
\begin{equation}\label{new sections}
\Sigma^\pm (b)=(\sqrt{\zeta_1 (b)}\cdot e^{i\theta_\pm
(b)},\sqrt{\zeta_1 (b)-b_2},\ldots ,\sqrt{\zeta_1 (b)-b_n}),
\end{equation}
are sections of $F$ which are smooth on $\R^n$. Let $\theta_0(b)=\Arg (ib_1)$. The section
\[
\Sigma^0 (b)=(\sqrt{\zeta_0 (b)}\cdot e^{i\theta_0(b)},\sqrt{\zeta_0 (b)-b_2},\ldots ,\sqrt{\zeta_0 (b)-b_n})
\]
is smooth on $\R^n\setminus\Delta$ and continuous on $\R^n$.
\end{cor}
\begin{proof}
It remains to verify that the above maps are sections of $F$. Let $\pi$ as in (\ref{c1 pi on V}). A direct computation shows that $\pi (\Sigma^\pm(b))=(\pm 1+ib_1,b_2,\cdots ,b_n)$ and $\pi (\Sigma^0(b))=(ib_1,b_2,\ldots ,b_n)$. Since $F$ factors via $\pi$ in an obvious way, the claim follows.
\end{proof}
Now let $\phi_i^t$ be the flow of the Hamiltonian vector field $V_{F_i}$ and consider the Poisson $\R^n$-action, $\Phi :\R^n\times\C^n\rightarrow \C^n$:
\begin{equation}\label{c1 eq action Phi}
\Phi(t_1,\ldots
,t_n;z)=\phi_1^{t_1}\circ\cdots\circ\phi_n^{t_n}(z).
\end{equation}

\begin{rem}\label{rm}
Observe that $\Phi$ is free and transitive along the fibres of $F$ over $\R^n\setminus\Delta$. Let $b_0\in\R^n\setminus\Delta$. Then, for each $z=\Sigma^-(b_0 )$ there is $(\alpha_1^0,\ldots ,\alpha_n^0)\in\R^n$ such that  $\Phi(\alpha_1^0,\ldots ,\alpha_n^0;z)\in\Sigma^+ (b_0)$. It follows from similar arguments to the ones used in the proof of Lemma
\ref{uniqueness of alpha} that there are locally defined
$C^\infty$ functions $\alpha_i(b)$ on $\R^n\setminus\Delta$ such that
$\alpha_i(b_0)=\alpha_i^0$ and such that $\Phi (\alpha_1(b),\ldots
,\alpha_n(b),z)\in\Sigma^+ (\R^n\setminus\Delta)$ for all $z\in \Sigma^-(\R^n\setminus\Delta)$.
\end{rem}

Denote by $\alpha :=\alpha_1$ and $\phi^t$ the flow of
$V_{F_1}$. Let $\mathcal O^-(b)$ and $\mathcal O^+(b)$ the
$T^{n-1}$-orbits of $\Sigma^-(b)$ and $\Sigma^+(b)$ respectively; it follows that  $\mathcal O^-(b)\cong\mathcal O^+(b)\cong T^{n-1}$.

\medskip
It is easy to see that for $x\in \mathcal O^-(b)$, $\phi^{\alpha
(b)}(x)\in\mathcal O^+(b)$. Let $z(b)\in\mathcal O^-(b)$ and
$w(b)\in\mathcal O^+(b)$. Let $\varphi$ denote the flow of
$\pi_\ast (V_{F_1})$. It is straightforward to check that the
solution to the equation $\varphi^{t(b)} (\pi (z(b)))=\pi (w(b))$
is precisely $t(b)=\alpha (b)$. We want to find an explicit
expression of $\alpha (b)$.
An easy computation shows that $\pi_\ast (V_{F_1})=-\chi\partial_u$ where,
\begin{equation}
\chi =\sum_{j=1}^n\frac{\prod_{i=1}^n |z_i|^2}{|z_j|^2}.
\end{equation}
Using formulae (\ref{circle-bundle-eq}), we see that $\chi
=\partial_xP_b(x)$. Regarding $\C\times\R^{n-1}=\R^{n+1}$ with
coordinates $(u,b)$, we can write $\pi_\ast (V_{F_1})$ as the
vector field in $\R^{n+1}$:
\begin{equation}
-2u\partial_xu\frac{\partial}{\partial u}
\end{equation}
Observe that for $b\in\R^n\setminus\Delta$ this vector field is
not singular.

\begin{lem}{\label{time-over-a-line}} Let $V$ be a vector field over $\R$. Let $p_0$ and $p$
be two points in $\R$ and assume $V(u)\neq 0$ for $u\in [p_0,p]$.
The time it takes to flow from $p_0$ to $p$ is:
\[T=\int_{p_0}^{p}\frac{du}{V(u)}.\]
\end{lem}

\begin{proof} Let $\varphi (t,u)$ be the flow of $V$. We want to find the time $T=T(p)$
such that $\varphi (T,p_0)=p$. We point out that
$\partial_t\varphi(t,p_0)|_{t=T(u)}=V(u)$. Then the derivative of $\varphi
(T(u),p_0)$ with respect to $u$ is $V(u)\partial_uT(u) =1$.
The claim follows easily from this.
\end{proof}

\begin{prop}{\label{alpha}} The function $\alpha$ is hypergeometric. Explicitly,
\begin{equation}\label{HL period}
\alpha(b)=-\int_{\zeta_0(b)}^{\zeta_1(b)}\frac{dx}{\sqrt{P_b(x)}},
\quad b\in\R^n\setminus\Delta,
\end{equation}
where $\zeta_0(b)$ is the maximal real root of
$P_b(x)=x(x-b_2)\cdots (x-b_n)-b_1^2$ and $\zeta_1(b)$ is the maximal real solution of $P_b(x)-1=0$.
\end{prop}

\begin{proof} First observe that $\pi (z(b))=(-1,b)$ and $\pi (w(b))=(1,b)$. It follows from Lemma \ref{time-over-a-line} that:
\begin{equation}\label{c1 ellip-int}
-\alpha
(b)=\int_{(-1,b)}^{(1,b)}\frac{du}{2u\partial_xu}.
\end{equation}
Bearing in mind that $u=\pm\sqrt{P_b(x)}$, it is not difficult to see
that $\alpha$ is as claimed.
\end{proof}
Of course, we can integrate $\alpha$ explicitly only when $n=2$. We are particularly interested in the case when $n\geq 3$, for which we need a precise understanding of the behaviour of $\alpha (b)$ as $b\rightarrow\Delta$. Let us write,
\begin{equation}
P_b=(x-\zeta_0)Q_b(x),\notag
\end{equation}
where $Q_b(x)$ is a polynomial with real coefficients. We notice
that $\zeta_0(b)$ becomes a (possible multiple) root of $Q_b(x)$
if and only if $b\in\Delta$.

\begin{prop}\label{c1 estimates} Let $\alpha$ as above and let
$\partial^k_{J_k}\alpha$ denote a partial derivative of order $k$.
Then, $\alpha$ is bounded from above by
\begin{equation}
-\frac{2}{(Q_b(\zeta_0))^{\frac{1}{2}}}.
\end{equation}
There are finite powers $w_0, w_1,\ldots , w_{n-1}\in\Z^+$,
depending on $J_k$, such that near $\Delta$,
\begin{equation}
|\partial^k_{J_k}\alpha |\simeq
\frac{1}{P'_b(\zeta_0)^{w_0}|\zeta_0-\beta_1|^{w_1}\cdots |\zeta_0-\beta_{n-1}|^{w_{n-1}}},
\end{equation}
where $\beta_i(b)={\rm Re}\, \rho_i(b)$ are the real part of the
roots of $Q_b$, $P_b=(x-\zeta_0)Q_b(x)$.
\end{prop}

\begin{proof} The proof involves the use of a (truncated) asymptotic expansion of $\alpha$. Since the integration limits of $\alpha$ depend on $b$, the estimates of $\partial^k_{J_k}\alpha$ turn out to be rather messy, as they involve the derivatives of $\zeta_0$. Here we estimate $\alpha$ to order $k=0$ and refer the reader to \cite{RCB} for the details concerning $k\geq 1$.

\medskip
Let $I=\alpha$ and let
$f=Q_b^{-\frac{1}{2}}$ and $dg=(x-\zeta_0)^{-\frac{1}{2}}dx$.
Integrating $I=\int fdg$ by parts we get,
\begin{equation}\label{k=0}
I=2\frac{(x-\zeta_0)^{\frac{1}{2}}}{Q_b(x)^{\frac{1}{2}}}\Bigg|_{\zeta_0}^{\zeta_1}-
\int_{\zeta_0}^{\zeta_1}-\frac{(x-\zeta_0)^{\frac{1}{2}}Q_b'(x)dx}{(Q_b(x))^{1+\frac{1}{2}}}.
\end{equation}

\medskip
Let $R_1$ be the first summand on (\ref{k=0}), and let $I_1$ be the
integral. We notice that $R_1=2(Q_b(\zeta_1))^{-\frac{1}{2}}$.
Since $x$ is such that $0\leq x-\zeta_0(b)\leq 1$, then
\begin{equation}
0\geq
I_1\geq\int_{\zeta_0}^{\zeta_1}-\frac{Q_b'(x)dx}{(Q_b(x))^{1+\frac{1}{2}}}
=2[(Q_b(\zeta_1))^{-\frac{1}{2}}-(Q_b(\zeta_0))^{-\frac{1}{2}}]
\end{equation}
Then we get, $|I|\simeq 2(Q_b(\zeta_0))^{-\frac{1}{2}}$.
\end{proof}

\begin{rem}\label{rem blowup of alpha} What
Proposition \ref{c1 estimates}  says is that the derivatives
of $\alpha$ blow up at $\Delta$ when the $\zeta_0$ becomes a
multiple root of $Q_b(x)$. Furthermore, $\alpha$ and all its
derivatives are bounded by a rational function having a pole of certain finite order along $\Delta$. For instance, when $n=3$, $\zeta_0$ becomes a root of $Q_b(x)$ as $b$ approaches to the spokes of $\Delta$, so $\alpha$ blows up there with order at most $-\frac{1}{2}$. This root becomes double at
$0\in\R^3$, so $\alpha$ blows up there with order at most $-1$.
\end{rem}

\medskip
Now let $B\subseteq\R^n$ be an open neighbourhood of $0\in\R^n$,
let $B_0:=B\setminus\Delta$ and let $\alpha_1,\ldots ,\alpha_n$ as in Remark \ref{rm}. Define a map $ A:\Sigma^-(B_0)\rightarrow\Sigma^+(B_0)$ as
$A(z)=\Phi (\alpha_1,\ldots ,\alpha_n; z)$.
In view of (\ref{new sections}), we can write $A$ explicitly as
$A:(z_1,z_2,\ldots ,z_n)\mapsto (-\bar z_1,z_2,\ldots z_n)$.
We verify that $A$ is smooth, furthermore, $A$ extends smoothly to $z\in\Sigma^-(B)$, regardless of the fact that the Poisson action is not freely transitive over singular fibres.

\medskip
Let $\tau_0=\sum\alpha_jdb_j$. We can find a 1-form $\eta=\sum\eta_jdb_j$ on $B$ such that $\tau:=\tau_0+\eta$ is closed. Indeed, let $\sigma$ be such that $d\sigma=\omega$ and let $\gamma (b)$ be a curve joining $\Sigma^-(b)$, and $\Sigma^+(b)$. One can verify that $\tau_0=dA_\gamma+R_\gamma$ where $A_\gamma=\int_\gamma\sigma$ and $R_\gamma$ is a 1-form  (c.f. (\ref{action int})). Defining $\eta=dH-R_\gamma$ for any $H\in C^\infty (B)$, we obtain $\tau=d(A_\gamma +H)$.

\medskip
Let $A':\Sigma^+(B)\rightarrow\Sigma (B) :=A' (\Sigma^+(B))$ be the map, $A'(z)=\Phi (\eta_1,\ldots ,\eta_n; z)$.
The composition $Q=A'\circ A$,
\begin{equation}\label{c1 map Q}
Q:\Sigma^-(B)\rightarrow \Sigma (B).
\end{equation}
is a $C^\infty$ map.
\medskip

\begin{prop}\label{ejemplo HL} Let $\tau=\tau_0+\eta$ be the 1-form as in the paragraph above. Then,
there is a symplectic manifold $(X, \omega)$ and a  Lagrangian
fibration $f:X\rightarrow B$ such that $\tau_1:=\tau$, and
$\tau_j=\pi db_j$, $j=1,\ldots n$, are the period 1-forms of $f$.
Furthermore, when $n=3$, $f$ coincides topologically with the
example in Construction \ref{c1 const 1,2}.
\end{prop}
\begin{proof}
We saw in Proposition \ref{c1 estimates} that the function $\alpha_1(b)$ is bounded from above by $-2(Q_b(\zeta_0(b)))^{-\frac{1}{2}}$.
We can find a smaller neighbourhood $B'\subseteq B$ of $\Delta$ such that $\alpha_1 (b) +\eta_1(b)<0$ for $b\in B'$. Let
$B'_0=B'\setminus\Delta$. Now let  $\mathcal O_b$ be the subset of $F^{-1}(b)$ defined by:
\[
\mathcal O_b:=\{\Phi (t;z)\mid t\in [\alpha_1 (b)+\eta_1(b),
0]\times [0,2\pi ]\cdots \times[0,2\pi ]\subseteq\R^n\} .
\]
Let $\overline{U}=\overline{\bigcup_{b\in B'}\mathcal O_b}$, this
is a $T^{n-1}$-invariant subset of $\C^n$. We see that for $b\in
B_0'$, $\overline{U}\cap F^{-1}(c)$ is a bounded cylinder and for
$b\in\Delta$ it is a bounded set in $F^{-1}(b)$. In both cases,
the boundary of these sets are the $T^{n-1}$-orbits: $\mathcal
T^-(b)=T\cdot\Sigma^- (b)$ and $\mathcal
T(b)=T\cdot\Sigma (b)$.

\medskip
Now let $W\subset F^{-1}(B')$ be a small $T^{n-1}$-invariant
neighbourhood of $\Sigma^-(B')$ such that $W\cap Crit
(F)=\emptyset$. Let $x\in W$ such that $F(x)=b$. There is a finite
$t\in\R^n$, $t=t(x)$, such that $x=\Phi (t,\Sigma^-(b))$.
Let $Q:\Sigma^-(B')\rightarrow \Sigma (B')$ as in
(\ref{c1 map Q}). Define a map $\mathcal Q:W\rightarrow
F^{-1}(B')$, $\mathcal Q(x)=y=\Phi (t(x),\Sigma (b))$. It
follows that $\mathcal Q$ extends $Q$. Moreover, similar arguments
to the ones used in Lemma \ref{uniqueness of alpha} can be used to
show that $t(x)$ is $C^\infty$. Let $W':=\mathcal Q (W)$. Then,
$\mathcal Q:W\rightarrow W'$ is clearly invertible, moreover,
$\mathcal Q$ is a diffeomorphism and $\mathcal Q$ sends $\mathcal
T^-$ diffeomorphically to $\mathcal T$.

\medskip
Now let $\mathcal U=W\cup\overline{U}\cup W'$ and let
$x,y\in\mathcal U$. We define $X=\mathcal U\slash\sim$ where
$x\sim y$ $\Leftrightarrow$ either $x=y$ or $y=\mathcal Q (x)$ if
$x\in W$ and $y\in W'$. By means of this identification, $X$ is a
smooth manifold. Intuitively, $\mathcal Q$ identifies the two
components on the boundary of $\bar U$. $F$ induces a smooth map $f:X\rightarrow
B'$ such that $f|_V=F|_V$ on a neighbourhood $V\subset int\; \bar
U$ of $Crit(f)$. If we can make $X$ into a symplectic manifold,
then the periods of $X$ are, by construction, $\tau_1,\ldots
,\tau_n$.

\medskip
Now let $\omega$ be the standard symplectic structure on $\C^n$.
Then $\omega$ restricts to symplectic forms on $W$ and $W'$. Let
us consider 1-form $\tau:=\tau_1$ which may be multi-valued on
$B'_0$. We can choose a domain $D\subset B_0'$ where $\tau$ is
single valued. Let $V_{\tau}$ the vector field determined by the
equation $F^\ast\tau=\iota (V_{\tau})\omega$. Since $\tau$ is
closed, $V_{\tau}$ is a symplectic vector field, i.e. its flow,
$\psi_s$, defines a 1-parameter family of symplectomorphisms. In
particular, its time $s=1$ flow map, $\psi_1$, is a
symplectomorphism. One can easily check that $\psi_1|_{W\cap
F^{-1}(D)}=\mathcal Q|_{W\cap F^{-1}(D)}$. This implies that
$\mathcal Q$ extends $\psi_1$ to $W$, in particular, $\mathcal
Q^\ast\omega$ and $\omega$ coincide on $W\cap F^{-1}(D)$. It
follows that $\mathcal Q^\ast\omega$ and $\omega$ coincide on
$W\setminus\{ f^{-1}(\Delta )\cap W\}\subset W$ which is dense in
$W$. Then, $\mathcal Q:W\rightarrow W'$ is a symplectomorphism and
therefore $X$ is a symplectic manifold. In dimension $n=3$, it is
not hard to see that $f$ coincides topologically with the example
given in Construction \ref{c1 const 1,2}.
\end{proof}

\begin{rem}\label{sil rem}
Observe that the symplectic structure in the example in Proposition
\ref{ejemplo HL} can be deformed by considering the 1-form
$\eta'=\eta +dH$ for any $C^\infty$ function $H$. It turns out that
there are cases for which $H$ does not induce a trivial
deformation of the symplectic structure. In fact, the example in
Proposition \ref{ejemplo HL} belongs to a large family of
Lagrangian fibrations whose members coincide topologically but
may not be symplectomorphic.
\end{rem}

\begin{thm}\label{c1 thm HL example} Let $\mathcal F=(X,\omega,f)\in \mathcal L(1,2)$. Then there are local sections $e_1,e_2,e_3$ of $R^2_cf^\#_\ast\Z$ such that the corresponding period 1-forms are:
\begin{align}
\tau_1 =\tau_0 +dH ,\quad \tau_2 = 2\pi db_2,\quad  \tau_3 =
2\pi db_3 \notag
\end{align}
where $\tau_0=\sum\alpha_idb_i$ is as in Proposition \ref{ejemplo HL},
$\alpha_1$ is as in (\ref{HL period}) and $H$ is a smooth function. Let $B\subseteq\R^3$ be an open ball. Secondly, for each $H\in C^\infty(B)$ there is a fibration $\mathcal F_H\in\mathcal L(1,2)$ with periods $\tau_1,\tau_2,\tau_3$ as above.
The monodromy representation of $\mathcal F\in\mathcal L(1,2)$ is
generated by the matrices: {\small $\begin{pmatrix}
  1 & 0 & 0\\
  1 & 1 & 0\\
  0 & 0 & 1
\end{pmatrix}$},
{\small $\begin{pmatrix}
  1 & 0 & 0\\
  0 & 1 & 0\\
  1 & 0 & 1
\end{pmatrix}$},
{\small $\begin{pmatrix}
  1 & 0 & 0\\
  1 & 1 & 0\\
  1 & 0 & 1
\end{pmatrix}$}.
\end{thm}
\begin{proof} The second statement follows from Proposition \ref{ejemplo HL} and Remark \ref{sil rem}. For the first claim, recall from Corollary  \ref{cor norm form 1,2} that any $\mathcal F\in\mathcal L(1,2)$ can be normalised in a neighbourhood $U\subset X$ of $p\in f^{-1}(0)\cap Crit(f)$ by $F:U\cong\C^3\rightarrow B\subset\R^3$, $F=(F_1,F_2,F_3)$ as in (\ref{c1 H-Lawson}). By redefining $X:=f^{-1}(f(U))$ if necessary, the restriction  $f|_{X\setminus U}$ induces a trivial bundle over $B$ with fibre $T^2\times [0,1]$. We can define sections $e_1,e_2,e_3\in R^2_cf^\#_\ast\Z$ in terms of the action of the Hamiltonian vector fields $v_i=v_{F_i}$ on $U$ and their extension to $X\setminus U$.  For $i=2,3$ and $b\in B$ we take $e_i(b)$ represented by integral curves $\gamma_i:[0,2\pi ]\rightarrow F^{-1}(b)$ of $v_i$. For $e_1(b)$, $b\in B\setminus\Delta$, we consider the sections $\Sigma_1:=\Sigma^+$ and $\Sigma_2:=\Sigma^-$ of $F$ as in Corollary \ref{cor sect HL} and define a representative $\gamma_1(b)$ of $e_1(b)$ as a suitable composition on flows of $v_1,v_2,v_3$, starting on $\Sigma_1(b)$ passing through $\Sigma_2(b)$ and returning to $\Sigma_1(b)$ in a completely analogous way as we did in the proof of Proposition \ref{prop. c1 period lattice}. The reader may easily check that the period 1-forms computed over $\gamma_i$ are $\tau_i$ as claimed.
\end{proof}

It is well known that the monodromy about the singular fibre of a focus-focus fibration can be explained in terms of a Dehn twist.  Similarly, for a fibration $\mathcal F\in\mathcal L(1,2)$, the monodromy of $\mathcal F$  can be understood as a ``two dimensional Dehn twist''.  For each generator of $\pi_1 (B\setminus\Delta,b)$, this twist is given by a full turn of a $T^2$-orbit, $\mathcal T (b)$, in one of the following ways: 
\begin{enumerate}
\item once in the direction of $v_2$; 
\item once in the direction of $v_3$; 
\item the turn in 1) followed by the turn in 2).
\end{enumerate}
\begin{figure}[!ht]
\begin{center}
\input{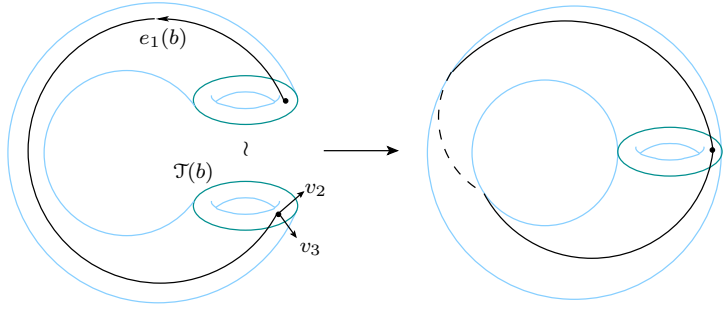}
\caption{Monodromy around a component of $\Delta -\{0\}$.}\label{Fig mon}
\end{center}
\end{figure}
In higher dimensions the description is analogous.

\section{The classification}
Let $\mathcal F=(X,\omega ,f)$ be a Lagrangian $T^3$ fibration over a smooth manifold $B$ and let $\Delta\subset B$ be the
discriminant locus of $f$. We are shall consider $\mathcal F\in\mathcal L(\kappa ):=\mathcal L(\kappa ,2)$, $\kappa=1,2$.

\begin{itemize}

\item Case $\kappa=1$: $B$ is an open 3-ball centred at $b_0\in\R^3$, $\Delta$ is a cone over 3 points.  Let $b_0\in B$ be the vertex of $\Delta$. There is only one singular point $p$ on the fibre $X_{b_0}$, i.e. the Poisson orbit $O_p=Crit(X_{b_0})=p$. There is a neighbourhood $U\subset X$ of $p$ and a normal form $\varphi\circ f|_U=q\circ\psi=F$ as in Corollary \ref{cor norm form 1,2}. The period lattice of $\mathcal F$ is as in Theorem \ref{c1 thm HL example}.

\item Case $\kappa=2$: $B=D\times (0,1)$, $\Delta=\{0\}\times (0,1)$. Let $b_0\in\Delta$ and $X_{b_0}$ the fibre over $b_0$. A point $p\in Crit (f)$ over $b_0$ belongs to a Poisson orbit $O_p=Crit (X_{b_0})\cong S^1$. There is a neighbourhood $U\subset X$ of $O_p$ and a normal form $\varphi\circ f|_U=q\circ\psi=F$ as in Definition \ref{defi one}. The period lattice of $\mathcal F$ is as in Proposition \ref{prop. c1 period lattice}.
\end{itemize}

\begin{defi} Let $(b_1,b_2,b_3)$ be coordinates on $B$ and let $\phi\in C^\infty (B)$. Let $\partial_{J_k}\phi$ denote a partial derivative of $\phi$. We say that $\phi$ is
\textit{$k$-flat} at $\Delta$ if $\partial_{J_k}\phi
|_b=0$ for each $b\in\Delta$ and each $J_k\leq k$. If $k=\infty$
we say that $\phi$ is \textit{flat} at $\Delta$.
Let $\varphi :B\subseteq\R^3\rightarrow \R^3$ be a $C^\infty$ map,
written as $\varphi =(\varphi_1,\varphi_2,\varphi_3)$. We say that
$\varphi$ is tangent to the identity at $\Delta$ if for each
$i=1,2,3$, the function $\varphi_i(b)-b_i$ is flat at $\Delta$.
\end{defi}

\begin{defi}\label{def ratio type} Let $\alpha$ be a function on $B$ which is $C^\infty$ on $B\setminus\Delta$.  We say that $\alpha$ is of
\textit{rational type} if for each $J_k\in\mathbb{Z}_{\geq 0}$ and any flat function $\phi$,
\[\lim_{b\rightarrow\Delta}\phi
\partial_{J_k}\alpha =0.\]
\end{defi}

\begin{ex} The function $\alpha (s,r)=\log |s|$ on $D\times (0,1)$ is of rational type. Similarly, $\alpha$ as in Proposition \ref{alpha} is also of rational type.
\end{ex}

\begin{defi} Let $\mathcal F, \mathcal F'\in \mathcal L(\kappa )$, $\kappa=1,2$, and let $\tau=\tau_0+dH$ and $\tau'=\tau_0+dH'$ be the singular periods of $\mathcal F$ and $\mathcal F'$ respectively. We say that $\mathcal F$ is \textit{formally equivalent} to $\mathcal F'$ if the function $H-H'$ is flat at $\Delta$.
\end{defi}

\begin{prop}\label{prop. c1 inv.} Let $\mathcal F=(X,\omega, f)$ and $\mathcal F'=(X',\omega',f')$ in $\mathcal
L(\kappa )$. Let $\tau_i$ and $\tau_i'$, $i=1,2,3$, be the
corresponding period 1-forms. If $\mathcal F$ and $\mathcal
F'$ are formally equivalent there is a
$C^{\infty}$-diffeomorphism between two small neighbourhoods of
$b_0\in\Delta$, $\varphi :B\rightarrow \varphi (B)=:B'$, such that
$\varphi^\ast (\tau_i')=\tau_i$ for $i=1,2,3$. Conversely, if there is a diffeomorphism $\varphi$ of $B$ matching $\tau_i$ and $\tau_i'$ and such that $\varphi$  is tangent to the identity, then $\mathcal F$ is formally equivalent to $(X',\omega',\varphi^{-1}\circ f')$.
\end{prop}
\begin{proof}  Let $\tau:=\tau_1$ and $\tau':=\tau_1'$ be the singular periods of $\mathcal F$ and $\mathcal F'$, expressed in the coordinates $(b_1,b_2,b_3)$ as in the previous section. We want to
find a diffeomorphism $\varphi$ such that
$\varphi^\ast\tau'=\tau$ and $\tau_j=\varphi^\ast\tau_j'$, $j=2,3$. The latter implies that $\varphi$ should be of the form:
\begin{equation}
b\mapsto (\varphi_1(b), b_2, b_3),
\end{equation}
where $\varphi_1$ is a smooth function to be determined.  Now, for $t\in[0,1]$, we
define a family of closed 1-forms $\tau_t=\tau +t(\tau'-\tau )$.
Then, $\tau_{t}|_{t=0}=\tau$ and $\tau_{t}|_{t=1}=\tau'$. Suppose
there is a 1-parameter family of maps, $G_t$, varying smoothly
with respect to $t$, such that each $G_t$ is a
diffeomorphism between small neighbourhoods of $0\in\C$
and such that $G_0$ is the identity map. Additionally, suppose
that:
\begin{equation}\label{eq. ff isotopy 1}
\frac{dG^\ast_t\tau_t}{dt}=0.
\end{equation}
Then, $G_1^\ast\tau_{t_1}=G_{0}^\ast\tau_{t_0}$ and we could
define $\varphi:=G_1$. It is standard to realise $G_t$ by means of integrating a time
dependent vector field $V_t$.
Using Cartan identity, we can rewrite equation (\ref{eq. ff
isotopy 1}) as:
\begin{align}\label{eq. ff
isotopy 1.1} G_t^\ast\left (\mathcal
L_{V_t}\tau_t+\frac{d\tau_t}{dt}\right ) = G^\ast_t\left
(d(\iota_{V_t}\tau_t)+\tau'-\tau\right ) =0.
\end{align}
Observe that $\tau' -\tau=d(H'-H)$. Then, the solution to
(\ref{eq. ff isotopy 1.1}) is determined by
\begin{equation}\label{eq. ff isotopy 2}
\iota_{V_t}\tau_t=H-H'.
\end{equation}
The solution should be of the form $V_t=g_t(b)\partial_{b_1}$
where $g_t(b)$ is a smooth function of $b$ and $t$. We observe that, since we want $\varphi
(\Delta )=\Delta$, then $V_t$ should satisfy $V_t(\Delta )=0$.
The left hand side of equation (\ref{eq. ff isotopy 2}) is:
\begin{align}
g_t(b)\cdot \Big(\alpha (b) +\frac{\partial}{\partial b_1} \big( H+t(H'-H)\big)\Big).
\end{align}
Define
\[
g_t(b)=\frac{H-H'}{\alpha (b) +\psi (b,t)}
\]
where $\psi(b,t)=\frac{\partial H}{\partial b_1}+t\frac{\partial
(H'-H)}{\partial b_1}$ is a $C^\infty$ function on $B\times [0,1]$. For $\kappa=1$ we know from Remark \ref{rem blowup of alpha} that $\alpha$ blows up at $\Delta$ with order at most $-1$. In the other hand, for $\kappa =2$ $\alpha$ blows up at $\Delta$ as a logarithm. Since $H-H'$ vanishes at $\Delta$ to all orders, for $\kappa=1,2$ we have 
\[
\lim_{b\rightarrow\Delta} \frac{(H-H')}{\alpha}=0.
\]
Therefore $g_t$ is continuous and $g_t(b)=0$ when $b\in\Delta$. In particular $V_t(\Delta)=0$ as required. A similar argument can be used to prove the smoothness of $g_t$. Indeed, for $\kappa=1$ the estimates carried out in Proposition \ref{c1 estimates} show that all the functions $\partial_{J_k}\alpha$ blow up with finite order along $\Delta$. This implies that for any  $h\in C^\infty(B)$ which is flat on $\Delta$,
\[
\lim_{b\rightarrow\Delta} h\partial_{J_k}\alpha =0.
\]
Now observe that the $k$-th partial derivatives of $(H-H')\slash\alpha$ are finite sums of terms of the type
\[
\frac{h\partial_{J_k}\alpha}{\alpha^m}
\]
with $m\leq 2k$ and $h\in C^\infty(B)$ flat at $\Delta$. It is not difficult to see from this that $\partial_{J_k}g_t$ is a continuous function on $B$ which vanishes at $\Delta$. A completely analogous argument is valid for the case $\kappa =2$.
Therefore $g_t(s)$ is a $C^\infty$ function on $B$ and it is flat at $\Delta$. This implies that the time one map of $V_t(b)$, $\varphi :=G_1$,
is a diffeomorphism which is tangent to the identity on $\Delta$.

\medskip
Now suppose there is a diffeomorphism $\varphi$ matching $\tau_i$ and $\tau_i'$ which is tangent to the identity at $\Delta$. Since $\varphi^\ast\tau_1'=\tau_1$ we can write $\varphi^\ast\tau_0-\tau_0=d(H-H'\circ\varphi )$. Furthermore, $\varphi$ can be written in coordinates $(b_1,b_2,b_3)$ as before as $\varphi=(\varphi_1,b_2,b_3)$, with $\varphi_1=\varphi_1(b)$ a smooth function on $B$. Observe that the 1-form $T:=\varphi^\ast\tau_0-\tau_0$ is single valued and smooth on $B$. Let us write $T=\sum T_idb_i$.  We claim that the functions $T_i$ are flat on $\Delta$. From this it follows directly that $H-H'\circ\varphi$ is flat at $\Delta$.  Now we will see that $\partial_{J_k}T_i|_\Delta=0$ for all $i=1,2,3$ and for each $J_k\leq k$, $k\in\mathbb{Z}_{\geq 0}$. Recall that $\tau_0=\sum\alpha_jdb_j$, where $\alpha_1=\alpha$ is a function on $B$ of rational type (cf. Definition \ref{def ratio type}) and $\alpha_2$ and $\alpha_3$ are locally defined. After an easy calculation we obtain
\begin{equation}\label{eq tees}
\begin{tabular}{l}
$T_1=(\alpha_1\circ\varphi)\partial_{b_1}\varphi_1-\alpha_1$\\
$T_2=(\alpha_1\circ\varphi)\partial_{b_2}\varphi_1 +\alpha_2\circ\varphi-\alpha_2$ \\
$T_3=(\alpha_1\circ\varphi)\partial_{b_3}\varphi_1 +\alpha_3\circ\varphi-\alpha_3$.
\end{tabular}
\end{equation}
Since $\varphi$ is tangent to the identity at $\Delta$, $\partial_{b_1}\varphi_1|_\Delta=1$ and $\partial_{b_1}\varphi_1(b)>0$ for $b$ in a small enough neighbourhood of $\Delta$. Then $T_1(b)\rightarrow 0$ as $b\rightarrow\Delta$ but since $T_1$ is continuous we have $T_1|_\Delta=0$. Similarly, $\partial_{b_2}\varphi_1$, $\partial_{b_3}\varphi_1$ and all the higher order derivatives of $\varphi_1$ vanish when restricted to $\Delta$. In particular $\partial_{b_2}\varphi_1$, $\partial_{b_3}\varphi_1$ are flat on $\Delta$. Then $(\alpha_1\circ\varphi)\partial_{b_j}\varphi_1|_\Delta=0$ for $j=2,3$. Now let us consider a (perhaps smaller) neighbourhood of $\Delta$ and a branch of $\alpha_j$ in this neighbourhood. Since $\varphi$ is infinitely tangent to the identity at $\Delta$, then $(\alpha_j\circ\varphi-\alpha_j)(b)\rightarrow 0$ as $b\rightarrow\Delta$. From the continuity of $T_j$ we have $T_j|_\Delta =0$ for $j=2,3$.  One can use the argument above inductively to show that $\partial_{J_k}T_i$ vanish on $\Delta$ for $k\geq 1$. 
\end{proof}

\begin{prop}\label{prop. c1 equiv} Let $\mathcal F$ and $\mathcal F'$ be in
$\mathcal L(\kappa )$. Let $\varphi :B\rightarrow \varphi
(B)=:B'$ be a diffeomorphism, such that $\varphi (\Delta)=\Delta$ and $\varphi^\ast (\tau_i')=\tau_i$ for $i=1,2,3$.
Then, there are sections $\Sigma$ and $\Sigma'$ of $f$ and $f'$
and a commutative diagram:
\[
\begin{CD}
X @>\Phi >> X' \\
@VVf V  @VVf' V \\
B @>\varphi>> B'
\end{CD}
\]
where  $\Phi$ is an orientation preserving diffeomorphism
sending $\Sigma$ to $\Sigma'$. The map $\Phi$ can be assumed to be
equivariant with respect to the $T^2$-actions induced by $\tau_j$
and $\tau_j'$ $j=2,3$. Furthermore, if $\Sigma$ and $\Sigma'$ are
Lagrangian, then $\Phi$ is a symplectomorphism.
\end{prop}
\begin{proof} Recall that $\mathcal F$ is normalised in a neighbourhood $U\subset X$ of the critical Poisson orbit $O_p\subset X_{b_0}$ by means of a symplectomorphism $\psi :U\rightarrow (V,\omega_0)$. Similarly, for $\mathcal F'$ there is a neighbourhood $U'\subset X'$ of $O_{p'}$ and a symplectomorphism $\psi':U'\rightarrow (V,\omega_0)$.

\medskip
Let $W\subseteq\psi (U)\cap\psi'(U')$. For simplicity
denote $U:=\psi^{-1}(W)$ and define $\Phi_0:=(\psi')^{-1}\circ\psi
|_U$.  Then, $\Phi_0$ is a symplectomorphism such that
$\Phi_0(O_p)=O_{p'}$ and such that $\varphi\circ F=F'\circ\Phi_0$.
Now let $\Sigma$ be a section of $F$ which does not pass through $O_p$.
Defining $\Sigma'=\Phi_0 (\Sigma )$ gives a section of $F'$ which
does not pass through $O_{p'}$. Notice that $\Sigma$ and $\Sigma'$
also define sections of $f$ and $f'$ respectively. Since $\Phi_0$
is a symplectomorphism, if $\Sigma$ is Lagrangian, then $\Sigma'$
is Lagrangian too.

\medskip
Let $\alpha$ be a local section of $T_{B}^\ast$. Let $v_{\alpha}$
be the vector field determined by the equation:
\begin{equation}\label{eq ff V_alpha}
F^\ast\alpha =\iota (v_{\alpha} )\omega .
\end{equation}
If we consider the 1-form, $db_i$, then $v_{db_i}=v_{q_i}$. As
we observed before, each $v_{q_i}$ extends to a vector field on $X$ which is tangent to the fibres of
$f$. Therefore, $v_{\alpha}$ extends to $X$ and, since the fibres
of $f$ are compact, the flow $g_{\alpha}^t$ of $v_{\alpha}$ is
defined for all $t\in\R$. For each $\alpha$ define the map
$T_{\alpha}:=g^{1}_{\alpha}:X\rightarrow X$. It follows that $\alpha\mapsto
T_\alpha$ induces a fibre preserving action, $T:T_{B}^\ast\times_B
X\rightarrow X$. Now define the map $\tilde\pi
:T_B^\ast\rightarrow X$ such that for each $\alpha_b\in
T_{B,b}^\ast$, $\tilde\pi (\alpha_b )=T_\alpha (\Sigma (b))=:x$,
which lies on the fibre $f^{-1}(b)$. One can verify that $x$ only depends
on the value of $\alpha$ at $b$. So, for $\bar\alpha$ such that $\bar\alpha
(b)=\alpha (b)$, $T_{\bar\alpha}(\Sigma (b))=T_\alpha (\Sigma
(b))=x$.

\medskip
Let $Z$ be the zero section on $T_B^\ast$. We know from Theorem
\ref{Thm. Jacobian} that $\tilde\pi (T_B^\ast )=X^\#$,
$\tilde\pi^{-1}(\tilde\pi(Z))=\Lambda$ is the period lattice of
$f$ and $\tilde\pi^{-1}|_{X^\#}:X^\#\rightarrow T_B^\ast$ is well
defined modulo $\Lambda$. Moreover, $\tilde\pi^{-1}|_{X^\#}$
composed with the projection $T_B^\ast\rightarrow
T^\ast_B\slash\Lambda =J_f$ gives a diffeomorphism $X^\#\cong
J_f$. If $\Sigma$ is Lagrangian this map is a symplectomorphism.

\medskip
Now let us take $\alpha'=(\varphi^{-1})^\ast\alpha$; this is a
local section of $T_{B'}^\ast$. Consider
the vector field $v_{\alpha'}$ induced by
$(F^{\prime})^\ast\alpha'=\iota (v_{\alpha'})\omega'$. Let
$g_{\alpha'}^t$ be the flow of $v_{\alpha'}$. Again, this flow is
complete, so we can define $T_{\alpha'}:X'\rightarrow X'$ such
that $T_{\alpha'}:=g_{\alpha'}^{1}$. Let
$\tilde\pi':T_{B'}^\ast\rightarrow X'$ such that $\tilde\pi'
(\alpha'_{b'})=T_{\alpha'}(\Sigma'(b') )=:x'$. Define the map
$\Phi^\# :X^\#\rightarrow X^{\prime\#}$ as the composition:
\begin{equation}
\xymatrix{
x\in X^\# \ar[rr]^{\Phi^\#} \ar[d]   & & x'\in X'^\# \\
[\alpha_b]\in J_f\ar[rr]_{(\varphi^{-1})^\ast} &  & [\alpha'_b]\in
J_{f'} \ar[u] }
\end{equation}
The horizontal map, which is induced by the pull back of sections
under $\varphi$, is well defined as $\varphi^\ast$ sends
$\Lambda'$ to $\Lambda$. The vertical maps, e.g. $[\alpha_b
]\mapsto g^{1}_\alpha (\Sigma (b))$, are independent of the choice
of the representative of $[\alpha_b ]\in J_f$. Indeed, let
$\tilde\alpha_b\in [\alpha_b ]$. Then, $\tilde\alpha_b =\alpha_b
+\lambda_b$, with $\lambda_b\in\Lambda_b$. It follows that
$g^t_{\alpha+\lambda}=g^t_\alpha$. In particular,
$g^1_{\alpha+\lambda}(\Sigma (b))=g^1_\alpha (\Sigma (b))$.

\medskip
We can write explicitly,
\begin{equation}\label{eq ff Phi sharp}
\Phi^\# (x)=g^{1}_{\alpha'}(\Sigma'(b'))
\end{equation}
where $x=g^{1}_{\alpha}(\Sigma (b))$ for some $[\alpha_b]\in J_f$
and $\alpha'=(\varphi^{-1})^\ast\alpha$. Notice that $\varphi$
induces a symplectomorphism between $T_{B'}^\ast$ and $T_B^\ast$.
Hence $\Phi^\#$ is a diffeomorphism and, when $\Sigma$ is
Lagrangian, $\Phi^\#$ is a symplectomorphism.

\medskip
Now let $X^\#\hookrightarrow X$ be the inclusion map and consider
$x\in U\cap X^\#$ over $b\in B$. We define
\[
  \Phi (x) =\begin{cases}
    \Phi^\#(x), & x\in X^\#, \\
    \Phi_0(x) & x\in U.
  \end{cases}
\]
The map $\Phi$ extends $\Phi^\#$ to $X$ and the
$T^2$-equivariance of $\Phi$ is verified \emph{a priori}. $\Phi$
is $C^\infty$ since the map $J_f\rightarrow J_{f'}$ is. We still need to check, however, that $\Phi^\# (x)=\Phi_0(x)$ on $U\cap X^\#$, i.e. that $\Phi$ is well defined. We prove this next.

\medskip
Let $x\in U\cap X^\#$ over $b\in B$ and define
$v_\alpha':= \Phi_{0\ast} (v_\alpha)$ and let $g^{t}_{v_\alpha'}$
denote the flow of $v_\alpha'$. We claim that the equation
$x=g^{1}_{\alpha}(\Sigma (b))$ implies that
\begin{equation}\label{eq ff yupi}
\Phi_0(x)=g^{1}_{v_\alpha'}(\Phi_0(\Sigma (b)))
\end{equation}
To see this let us regard $\gamma(t):=g_\alpha^t(\Sigma (b))$ as
the integral curve of $v_\alpha$ such that $\gamma (0)=\Sigma (b)$
and $\gamma (1)=x$. Now let $\gamma' (t):=\Phi_0 \left (\gamma
(t)\right )$. This is a curve on $F^{\prime -1}(b')$, $b'=\varphi
(b)$, such that $\gamma'(0)=\Phi_0(\Sigma(b))=\Sigma'(b')$ and
$\gamma'(1)=\Phi_0(x)$. Furthermore, $\gamma'$ is an integral
curve of $v_\alpha'$. Indeed, we see that:
\[
\frac{d\gamma'}{dt}=\frac{d(\Phi_{0}\circ\gamma )}{dt}=
\Phi_{0\ast}( v_\alpha )=v_\alpha'.
\]
Therefore $\gamma'(t)=g^t_{v_\alpha'}$ and $g^1_{v_\alpha'}\left (
\Sigma'(b')\right )=\Phi_0(x)$. Now observe that
$\Phi_0^\ast\omega'=\omega$ implies that:
\begin{equation}\label{eq ff v v'}
v_\alpha'=v_{\alpha'}.
\end{equation}
To prove this we notice that $F^\ast\alpha
=\Phi_0^\ast(F'^\ast\alpha')$. Now we can write (\ref{eq ff
V_alpha}) as:
\begin{align}\label{ff eq dhr}
\Phi_0^{-1\ast}\left (\iota (v_\alpha )\Phi_0^\ast\omega'\right
)=F^{\prime\ast}\alpha'
\end{align}
The left hand side of (\ref{ff eq dhr}) can be written as $\iota
(\Phi_{0\ast}v_\alpha )\omega'$. Then, it follows that
$v_\alpha'=\Phi_{0\ast}v_\alpha =v_{\alpha'}$. Now, from  (\ref{eq
ff yupi}) and (\ref{eq ff v v'}) we conclude that:
\[
\Phi_0(x)=g^1_{\alpha'}(\Sigma'(b'))
\]
which is equal to $\Phi^\# (x)$ in (\ref{eq ff Phi sharp}), hence $\Phi$ is well defined.

\medskip
Observe that, for $\kappa=2$, we can start the above construction in terms of the section $\Sigma=\Sigma_1$ as in Construction \ref{contr. action}, which is Lagrangian. Therefore $\Phi$ turns out to be a symplectomorphism. In the case $\kappa=1$, the sections $\Sigma^\pm$ as in (\ref{new sections}) are not Lagrangian. This does not give much trouble as we can always find a Lagrangian section. The argument is valid for $\kappa =1,2$. Let $U\subset X$ be as before. Observe that for any given section $\Sigma_0$ of $f$ with $\Sigma_0(B)\subset U$, there exists a neighbourhood $\mathcal U\subseteq U\subset X$ of $\Sigma_0$ such that $\mathcal U\cap Crit(f)=\varnothing$ and a fibre-preserving symplectomorphism $(\mathcal U,\omega|_\mathcal{U})\rightarrow(T_B^\ast,\Omega )$. Here $\Omega =\omega_0+\beta$ where $\omega_0$ is the standard symplectic structure on $T_B^\ast$ and $\beta$ is the pull-back under $T_B^\ast\rightarrow B$ of a closed 2-form on $B$ (if $\Sigma_0$ were Lagrangian $\beta=0$). Observe that in our situation $\beta$ can be assumed to be exact, so we have $d\theta=\Omega -\omega_0$ for some 1-form $\theta$ on $B$. Then $-\theta$ defines a section, $\Sigma_\theta$, of $T_B^\ast$ which is Lagrangian with respect to $\Omega$. Then $\Sigma_\theta$ maps to a Lagrangian section, $\Sigma$, of $f$ inside $\mathcal U$. Using $\Sigma$ to define $\Phi$ we obtain a symplectomorphism.
\end{proof}

\begin{thm} Let $\mathcal F =(X,\omega, f)$ and $\mathcal F'=(X',\omega',f')$ be Lagrangian fibrations of type $\mathcal L(\kappa )$, $\kappa =1,2$. Then $\mathcal F$ is formally equivalent to $\mathcal F'$ if and only if $\mathcal F$ is symplectically equivalent to $\mathcal F'$.
\end{thm}
\begin{proof} Assume $\mathcal F$ and $\mathcal F'$ are formally equivalent. Then Proposition \ref{prop. c1 inv.} gives us a diffeomorphism $\varphi$ on $B$ such that $\varphi^\ast\tau_j'=\tau_j$. In view of Proposition \ref{prop. c1 equiv}, $\varphi$ lifts to a fibre-preserving symplectomorphism $\Phi :X\rightarrow X'$.

\medskip
To prove the converse we suppose there is a symplectomorphism $\Psi$ and a suitable diffeomorphism $\varphi$, making a commutative diagram:
\begin{equation}\label{eq ff last diag}
\xymatrix{
  X \ar[rr]^{\Psi} \ar[d]_{f}
                &  &    X' \ar[d]^{f'}    \\
                B \ar[rr]^{\varphi} & &B'=\varphi (B) }
\end{equation}
One can always take a diffeomorphism, $\tilde\varphi$,
from a neighbourhood of $b_0\in B$ into a neighbourhood, $\tilde
B\subseteq B'$, of $\varphi (b_0)=b_0'$ and such that
$\varphi\circ\tilde\varphi ^{-1}$ is tangent to the
identity at $\Delta\cap\tilde B$. Let $\tilde f=
\tilde\varphi\circ f$. Then, $(X,\omega ,\tilde f)$ and $\mathcal F$ define the same germ. Now, $\Psi$ together with the map
$\varphi':=\varphi\circ\tilde\varphi^{-1}:\tilde B\rightarrow
\varphi'(\tilde B)\subseteq B'$ makes $(X,\omega ,\tilde f)$ and
$\mathcal F'$ symplectomorphic, with $\varphi'$ being
tangent to the identity at $\Delta$. Let us denote $f:=\tilde f$
and $\varphi:=\varphi'$.

\medskip
We claim now that $\tau_i=\varphi^\ast\tau_i'$. To see this we
take $V_i'$ to be the vector fields determined by the equation
\begin{equation}\label{eq. ff periods again}
f^{\prime\ast}\tau_i'=\iota_{V_i'}\omega',\qquad i=1,2,3 .
\end{equation}
These vector fields are defined on open sets $f^{\prime -1}(U')$,
where $U'\subset B'_0$ is an open set on which a branch of
$\tau_1'$ is defined. It follows that $V_i'$ are vector
fields whose flows are periodic. We can take integral curves of
$V_i'$ to define simple loops, $\gamma_i'(b)$ and
on $f^{\prime -1}(b)$, representing the cycles
$e_i'(b)$ generating $H_1(f^{\prime -1}(b'),\R)$. These loops can be used for computing the
period 1-forms of $f'$ which are, tautologically, $\tau_i'$. Now
define $V_i$ to be the vector fields determined by the equation
$\iota (V_i)\omega=(f'\circ\Psi )^\ast\tau_i$. Since $\Psi$ is
symplectic, $\Psi_\ast V_i=V_i'$. The above implies that the flow of $V_i$ is
periodic. One verifies that suitable integral curves $\gamma_i$ of
$V_i$ generate $H_1(X_b,\Z )$, so we can define the period one
forms, $\tau_i$  of $f$ by integrating along $\gamma_i$. Now
observe that, since the diagram (\ref{eq ff last diag}) commutes,
$V_i$ also satisfies the equation $\iota (V_i)\omega
=\varphi^\ast\tau_i'$. Therefore, $\tau_i=\varphi^\ast\tau_i'$.
The conclusion follows now from Proposition \ref{prop. c1 inv.}.
\end{proof}

\subsection{Acknowledgements} I would like to thank M. Gross, my thesis advisor, for his help and support. I am also grateful to D. Matessi, J. Rawnsley and R. Thomas for their useful comments and suggestions. I would like to thank the Mathematics Institute of Warwick University, U.K. and the I.C.T.P. in Trieste.
I also wish to thank the referee of this article for making helpful suggestions that improved the final version of this paper, specially for pointing out to me a way to simplify \S 5.
\bibliographystyle{plain}
\bibliography{/home/rcastano/TEX/biblio}

\bigskip
\begin{flushleft}
Address:\\
Ricardo Casta\~no-Bernard \\
International Centre for Theoretical Physics\\
Mathematics Section\\
Strada Costiera 11\\
Trieste 34014, Italy

\medskip
e-mail: rcastano@ictp.trieste.it
\end{flushleft}

\end{document}